\documentclass[11pt,reqno,oneside]{amsart}
\usepackage{verbatim,amsmath,amssymb,cite,epsfig,xspace,color,subfigure,bbm}
\usepackage[foot]{amsaddr}

\numberwithin{equation}{section}

\theoremstyle{plain}
\newtheorem{theorem}{Theorem}[section]

\theoremstyle{definition}

\theoremstyle{remark}
\newtheorem{remark}{Remark}[section]

\newcommand{\calL}{\mathcal{L}}
\newcommand{\calF}{\mathcal{F}}

\newcommand{\calI}{\mathcal{I}}
\newcommand{\calP}{\mathcal{P}}
\newcommand{\calQ}{\mathcal{Q}}
\newcommand{\calS}{\mathcal{S}}
\newcommand{\calM}{\mathcal{M}}

\newcommand{\mbf}[1]{\mathbf{#1}}
\newcommand{\mbb}[1]{\mathbb{#1}}

\begin{document}

\title[Coarse--graining of overdamped Langevin dynamics]{Coarse-graining of overdamped Langevin dynamics via the Mori-Zwanzig formalism}

\author{Thomas Hudson}
\address{Mathematics Institute, University of Warwick, Coventry, CV4 7AL}
\thanks{The work of T. Hudson is supported by the Leverhulme Trust through Early Career Fellowship ECF-2016-526.}

\author{Xingjie Helen Li}
\address{Fretwell 350C, 9201 Univ City Blvd., Charlotte, NC, 28023}
\thanks{The work of X. Li is supported in
  part by the Simons Collaboration Grant with Award ID: 426935 and NSF DMS-1720245.}

\date{}
\begin{abstract}
  The Mori--Zwanzig formalism is applied to derive an equation for the evolution of linear
  observables of the overdamped Langevin equation. To illustrate the resulting equation and its
  use in deriving approximate models, a particular benchmark example is studied both numerically
  and via a formal asymptotic expansion. The example considered demonstrates the importance
  of memory effects in determining the correct temporal behaviour of such systems.
\end{abstract}

\maketitle

\section{Introduction}
Molecular dynamics (MD) is a widely--used simulation technique which captures the atomistic details
of material systems, allowing the prediction of their properties and behavior
\cite{Mazenko2006,Evans2008}. However, despite the vast increases in computational capacity over recent decades, it is still not always possible to work with MD models at
full resolution, particularly when studying large, complex systems over long time--scales.
Fortunately, in many cases, the objectives of a simulation occur within a small region of interest.
This observation has led to the development of \emph{coarse--grained MD} (CGMD) models, in which
excess degrees of freedom are incorporated implicitly
\cite{Ford1965a,Munakata1985a,Espanol1995a,Majda2005a,Izvekov2006a,Kinjo2007a,Snook2007a,
  Engquist2007,Ebook2011a}.

Building reliable and efficient CGMD models attuned to the quantities of interest is a difficult
problem. First, the simulator must find appropriate variables which capture the quantities of
interest \cite{Guttenberg2013a}, often termed \emph{reaction coordinates}. Once these are fixed,
an appropriate proxy model for the reaction coordinates must be obtained, which implicitly
incorporates the interaction between the reaction coordinates and unresolved degrees of freedom
\cite{Mazenko2006,Engquist2007,Lelievre2010,Ebook2011a,Legoll2012,Velinova2017}.
If the objectives of a simulation are `static' macroscopic equilibrium properties such as free
energy or reaction rates, then a wide variety of choices of proxy dynamics which appropriately
sample the relevant measures are available. However, if the objective is to capture a dynamical
property of the physical system such as kinematic viscosity or a diffusion rate, then it is
important to capture the correct effective dynamics of the reaction coordinates arising due to
the relevant dynamics of the full system over moderate timescales \cite{Majda2006a,Mazenko2006,
  Hartmann2007,Evans2008,Chorin2013a,Swinburne2015a}.

In recent years, a variety of studies of CGMD schemes have been undertaken, aiming to analyse the
predictions of such schemes. In all cases, the ultimate goal is to obtain verifiable, statistically
accurate predictions of the true dynamics for various applications. The wide variety of
mathematical techniques used includes
\begin{itemize}
\item optimal prediction techniques \cite{Chorin2000a,Chorin2002a,Chorin2006a,Frank2011a};
\item information-theoretic tools \cite{Majda2005a,Majda2012a,Plechac2013a,Majda2016a,Dobson2016a};
\item statistical filtering and ensemble methods \cite{Majda2004a,Majda2015a,FeiLu2017a};
\item identification of an appropriate parametrization \cite{Lei2010a,Majda2011a,Majda2013a,FeiLu2015a};
\item series expansion \cite{Pinski2014a,Venturi2017a,Lei2016a};
\item pathwise estimates \cite{Legoll2017a,Lelievre2018a,Klus2018,Legoll2018b}; and
\item conditional expectations \cite{Legoll2010a}.
\end{itemize}
Here, our focus is on the Mori--Zwanzig (MZ) approach to CGMD benchmark problem
\cite{Mori1965,Zwanzig1973a,Berkowitz1983a,Berne2000a,Zwanzig2001a}. The MZ formalism provides
an exact expression of the dynamics for a CGMD scheme, and is governed by three
terms which separate out different contributions to the true dynamics, each of which has
a different statistical physical meaning. This decomposition allows a study of the
sources of error: the first term accounts for a conservative dynamics
due to the effective interactions between the coarse grained variables;
the second is a history-dependent term determined by a time integral of a memory
kernel which represents the interactions between the resolved and unresolved variables; and the
third term represents the random thermal fluctuations arising from unresolved variables.
In different situations, each of these terms may have a more or less important role, but to
correctly capture the dynamical properties and validate an effective model, it is critical to
measure the relative size and behaviour of these terms accurately.

Our study concentrates particularly on the \emph{memory term}, which may heuristically be thought of
as measuring the extent to which the set of reaction coordinates is decoupled from the unresolved
degrees of freedom. In recent years, there have been tremendous efforts to investigate
memory terms from MZ projections, for a variety of classes of dynamics, see for example
\cite{Berkowitz1983a,Chorin2000a,Darve2009a,XiantaoLi2010a,Yoshimoto2013a,Guttenberg2013a,
  XiantaoLi2014a,ZhenLi2015a,Stinis2015a,Venturi2016a,ZhenLi2017a,Venturi2017b}.
One common approach is to hope that a time-scale separation between the resolved and unresolved
variables occurs, i.e. the fluctuations of the unresolved variables occur on a much faster
timescale than those of the resolved variables, and therefore the two sets of variables are weakly
correlated. In such cases, the memory kernel decay rapidly, approximating a delta function in
time \cite{Chorin2013a}.

Our aim in this paper is to demonstrate that while such delta approximations of the memory kernel
are appropriate in many situations, it is not generally to be expected that the
memory kernel is independent of the value of the reaction coordinate, even in the simple situation
where the chosen reaction coordinates are linear. To capture the
correct dynamics, further careful analysis and sampling of the memory is therefore required.

As
an illustration of
this issue, we consider the dynamical behaviour of a gradient flow with stochastic forcing (often
called the overdamped Langevin equation), demonstrating that at least in this case,
a na\"ive approach to approximating the memory kernel yields a poor approximation of the dynamics.
We hope that the benchmark problem we consider here will provide insight which will enable the study of CGMD derived from full Langevin dynamics based on reliable asymptotic analysis
in future.

\subsection{Outline}
This paper is organized as follows. In Section~\ref{sec:prob_formulation}, we review the
Mori--Zwanzig formalism applied to general gradient flow systems, and in
Theorem~\ref{thm:MZ_formulation}, derive an exact equation for the evolution of linear observables
within an abstract framework. Our benchmark example is discussed in Section~\ref{sec:benchmark},
and an asymptotic analysis is performed to obtain approximations of the various terms in the MZ
equation. Finally, we study this particular example numerically in Section~\ref{sec:numerics}.

\section{Formulation of the problem}\label{sec:prob_formulation}
As our reference fine--scale dynamical system, we consider the following overdamped Langevin dynamics defined on $\mbb{R}^N$:
\begin{equation}\label{over_damped_Langevin}
d\mbf{X}_t=
-\nabla_{\mbf{x}} V(\mbf{X}_t) dt +\sqrt{2\beta^{-1}}d \mbf{B}_t.
\end{equation}
Here, $\mbf{B}_t$ denotes a standard vector--valued Brownian motion, $V(\cdot)$ is a potential
energy and $\beta$ is the inverse temperature. Throughout this work, we assume that $V$ is at least
of class $C^2$, and satisfies
the following conditions:
\begin{enumerate}
\item There exist constants $\gamma>0$ and $\delta\geq 0$ and $R>0$ such that
  \begin{equation*}
    -\nabla V(\mbf{x})\cdot \mbf{x} \leq \gamma|\mbf{x}|^2+\delta \qquad\text{for all }|\mbf{x}|\geq R.
  \end{equation*}
\item The gradient $\nabla V(\mbf{x})$ is globally Lipschitz, i.e. there exists $\alpha>0$ such
  that
  \begin{equation*}
    \big|\nabla V(\mbf{x})-\nabla V(\mbf{y})\big|\leq \alpha\,|\mbf{x}-\mbf{y}|.
  \end{equation*}
\end{enumerate}
Under condition (1), it is well--known \cite{RT96} that the dynamics defined by
\eqref{over_damped_Langevin} are ergodic with respect to the Gibbs measure, $\mu$, given by
$d\mu(\mbf{x})=\frac1Ze^{-\beta V(\mbf{x})}d\mbf{x}$ where $Z$ is the partition function
\begin{equation*}
  Z:=\int e^{-\beta V(\mbf{x})}d\mbf{x}.
\end{equation*}

Given a regular function $\mbf{F}:
\mbb{R}^N \rightarrow \mbb{R}^m$, which may be thought of as describing a family of \emph{reaction coordinates}, we may apply It\^o's formula to deduce that the value of
$\mbf{F}(\mbf{X}_t)$ is governed by the It\^o SDE
\begin{equation}\label{Liouville_eq1}
d\mbf{F}(\mbf{X}_t) = \Big(-\nabla V(\mbf{X}_t)\cdot \nabla \mbf{F}(\mbf{X}_t)+\beta^{-1}\Delta \mbf{F}(\mbf{X}_t)  \Big) dt
+\sqrt{2\beta^{-1}}\nabla \mbf{F}(\mbf{X}_t)\cdot d\mbf{B}_t.
\end{equation}
In particular, if we consider a linear coarse--graining selector $\mbf{F}(\mbf{x}):= \Phi \mbf{x}$, where $\Phi \in \mbb{R}^{m\times N}$ is a constant matrix, \eqref{Liouville_eq1} becomes
\begin{equation}\label{Liouville_eq2}
d\mbf{F}(\mbf{X}_t) = -\Phi \nabla V(\mbf{X}_t)  dt
+\sqrt{2\beta^{-1}}\Phi \, d\mbf{B}_t.
\end{equation}
Such linear coordinates are commonly used in CGMD schemes, particularly for large molecules such
as polymers \cite{Espanol1995a,dPMC12,Pasquale2018}.
If we are interested in the value of the reaction coordinates described by $\mbf{F}$ alone, then
\eqref{Liouville_eq2} provides an equation for their evolution. In general however, since the
first term on the right--hand side of the equation depends on the full process $\mbf{X}_t$, this
is not a closed equation for the value of $\mbf{F}(\mbf{X}_t)$.

In order to formulate a closed approximate equation for $\mbf{F}(\mbf{X}_t)$, we use the
\emph{Mori--Zwanzig formalism}, which uses projection operators to decompose the equations
describing observables of a dynamical system into terms involving the value of the observables
alone, and `error' terms describing the contribution of variations of $\mbf{X}_t$ which do not
directly change the value of the observable.

In this case, the projection operator we choose to apply is the \emph{Zwanzig projection}, which
involves taking a conditional expectation with respect to the Gibbs distribution, i.e.
\begin{equation}\label{MZ_proj_def}
\calP \mbf{g}= \mbb{E}_\mu\big[\mbf{g}\,\big|\, \mbf{F}(\mbf{x})= \mbf{h} \big]
: = \frac {\int_{\mbf{F}^{-1}(\mbf{h})}\mbf{g}(\mbf{x})\, e^{-\beta V(\mbf{x})} d\mbf{x}  }{\int_{\mbf{F}^{-1}(\mbf{h})} e^{-\beta V(\mbf{x})} d\mbf{x}};
\end{equation}
note that in the above formula, we have cancelled the common factor $\frac1Z$ from the numerator
and denominator.\footnote{Note that to be completely technically correct, the right hand side should
  be understood in the sense of Radon--Nikodym differentiation of measures.} It may be verified that $\calP$ is an orthogonal projection on the space of square--integrable observables, i.e. $L^2(\mbb{R}^N;e^{-\beta V(\mbf{x})}d\mbf{x})$, and we can therefore define its orthogonal counterpart,
\[
\calQ:= \calI-\calP.
\]

In particular, we note that the evolution of \eqref{Liouville_eq2} can be divided into the
stationary, mean--zero process induced by the Brownian motion, and the evolution of the mean
value of $\mbf{F}(\mbf{X}_t)$. To consider the behaviour of the latter quantity given knowledge
of $\mbf{X}_0$, we define
\begin{equation}\label{observable}
  \mbf{h}_t(\mbf{x})=\mbb{E}\big[\mbf{F}(\mbf{X}_t)\,\big|\,\mbf{X}_0 = \mbf{x}].
\end{equation}
The evolution of this quantity is governed by the usual generator of the SDE
\eqref{over_damped_Langevin},
\begin{equation}\label{Liouville_def_gradient}
\calL:= -\nabla V\cdot \nabla+\beta^{-1}\Delta.
\end{equation}
Using this definition, the Feynmann--Kac formula governing the evolution of $\mbf{h}$ states that
the function $\mbf{h}$ solves the PDE
\begin{equation}\label{eq:FK}
  \partial_t\mbf{h}_t=\calL\mbf{h}_t\qquad\text{with}\qquad\mbf{h}_0=\mbf{F}.
\end{equation}
Using the definition of $\calP$, we apply the Mori--Zwanzig formalism to provide a different
expression of \eqref{eq:FK}, stated in the following theorem.

\begin{theorem}\label{thm:MZ_formulation}
Let $\mbf{X}_t$ satisfy the SDE on $\mbb{R}^N$
\[
d\mbf{X}_t=-\nabla V(\mbf{X}_t)dt+\sqrt{2\beta^{-1}}d\mbf{B}_t,
\]
and given a constant matrix of full rank $\Phi\in \mbb{R}^{N\times m}$, the observable
\begin{equation}\label{observable}
  \mbf{h}_t(\mbf{x})=\mbb{E}\big[\mbf{F}(\mbf{X}_t)\,\big|\,\mbf{X}_0 = \mbf{x}];
\end{equation}
satisfies the following integro--differential equation:
\begin{equation}\label{MZ_formulation}
  \partial_t \mbf{h}_t(\mbf{x})=-\Phi\Phi^T\,\nabla \calS\big(\mbf{h}_t(\mbf{x})\big)
  +\int_{0}^{t} \calM_s\big(\mbf{h}_{t-s}(\mbf{x})\big)\cdot \nabla \calS\big(\mbf{h}_{t-s}(\mbf{x})\big)-\frac{1}{\beta}\mathsf{div}\,\calM_s \big(\mbf{h}_{t-s}(\mbf{x})\big)ds+\calF_t(\mbf{x}),
\end{equation}
where:
\begin{enumerate}
\item $\calS:\mbb{R}^m\to\mbb{R}$ is the \emph{effective potential}, defined to be
  \begin{equation}\label{effective_potential}
    \calS(\mbf{h}):= -\frac{1}{\beta} \log Z_\Phi(\mbf{h})
    \text{ with }Z_\Phi(\mbf{h}):=\int_{\mbf{F}^{-1}(\mbf{h})} e^{-\beta V(\mbf{x})} d\mbf{x},
  \end{equation}
\item $\calM_s:\mbb{R}^m\to\mbb{R}^{m\times m}$ is the \emph{memory kernel}, defined to be
  \begin{equation}\label{memory_eq2}
\calM_s(\mbf{h}):=\beta \,\mbb{E}\big[e^{s\calQ\calL}\calQ\calL \mbf{F} \otimes \calQ\calL \mbf{F}\,\big|\, \mbf{F}(\mbf{x})=\mbf{h}\big]=\beta \,\mbb{E}\big[\calF_s \otimes \calF_0\,\big|\, \mbf{F}(\mbf{x})=\mbf{h}\big],
\end{equation}
\item and $\calF_t:\mbb{R}^N\to\mbb{R}$ is the \emph{fluctuating force}, defined to be
  \begin{equation}\label{RandomForce}
    \calF_t :=  e^{t\calQ\calL}\calQ\calL\mbf{F}.
  \end{equation}
\end{enumerate}
\end{theorem}
\medskip
A proof of this result is given in Appendix~\ref{sec:Proof}, and involves adapting standard variants of
the Mori--Zwanzig formalism already present in the literature to this stochastic setting.

\begin{remark}
  The Mori--Zwanzig formalism \cite{Zwanzig1961,Mori1965,Zwanzig1973a,Nordholm1975} uses projection
  operators to rewrite the
  equations governing observables of a dynamical system. Various formulations have been developed
  in recent years with a variety of applications in mind, and influence our own derivation,
  including: \cite{XiantaoLi2010a} treating crystalline solids via the harmonic approximation;
  \cite{Venturi2017a} for the harmonic oscillators based on operator series expansions of the
  orthogonal dynamics propagator; \cite{Ma2018a} for the full Langevin dynamics model based on
  reduced-order modeling; \cite{Hijon2010a} for a model based on dissipative particle dynamics;
  and \cite{Pasquale2018} for a `hybrid' coarse--graining map of a Hamiltonian model.

  Recombining the evolution of the mean $\mbf{h}_t$ given by \eqref{MZ_formulation} and adding back
  the Brownian motion, we find that \eqref{Liouville_eq2} can be written
  \begin{equation}\label{MZStochastic}
    \begin{aligned}
    d\mbf{F}(\mbf{X}_t) &= -\Phi\Phi^T\,\nabla \calS\big(\mbf{F}(\mbf{X}_t)\big)dt
    +\int_{0}^{t} \calM_s\big(\mbf{F}(\mbf{X}_{t-s})\big)\cdot \nabla \calS\big(\mbf{F}(\mbf{X}_{t-s})\big)ds\,dt\\&\hspace{2cm}-\int_0^t\frac{1}{\beta}\mathsf{div}\,\calM_s \big(\mbf{F}(\mbf{X}_{t-s})\big)ds\,dt
    +d\calF_t+\sqrt{2\beta^{-1}}\Phi d\mbf{B}_t.
    \end{aligned}
  \end{equation}
  It is important to note at this point that \eqref{MZStochastic} is equivalent
  to considering the full evolution \eqref{Liouville_eq2}, in particular because $\calF_t$ is
  unknown. Equation \eqref{MZStochastic} therefore remains
  unclosed; however, the power of this formulation is that if $\calF_t$ has statistics which are
  well--captured by some proxy process $\widetilde{\calF}_t$, and $\calS$ is either known or
  accurately approximated, then we can obtain a closed--form approximate dynamics
\begin{multline*}
  d\mbf{h}_t = -\Phi\Phi^T\nabla\calS(\mbf{h}_t)dt +\bigg(\int_0^t\widetilde{\calM}_s(\mbf{h}_{t-s})\nabla \calS(\mbf{h}_{t-s})-\frac1\beta\mathsf{div}\,\widetilde{\calM}_s(\mbf{h}_{t-s})ds\bigg)\,dt\\
  +d\widetilde{\calF}_t +\sqrt{2\beta^{-1}}\Phi d\mbf{B}_t,
\end{multline*}
where $\widetilde{\calM}_s(\mbf{h})$ is the autocovariance function of  $\widetilde{\calF}_t(\mbf{h})$ (see for example \cite{Lindgren2012}).
\end{remark}

To explore this formulation and better understand the relationship between the terms involved, the
remainder of the paper is devoted to an exploration of a dynamics where an accurate approximation
\eqref{MZ_formulation} and the terms within
it.

\section{A benchmark problem}\label{sec:benchmark}
In this section, we will consider the overdamped Langevin equation \eqref{over_damped_Langevin}
in the particular case where $\mbf{x}\in\mbb{R}^2$, and the potential energy is defined to be
\begin{equation}\label{potential_eq}
  V(x,y):=\frac{\mu}{2}x^2+\frac{\lambda}{2}\big(\tau \sin(\omega x)-y\big)^2
\end{equation}
with $\mu$, $\lambda$, $\tau$, $\omega\geq0$ being parameters. Specifying even further, we will
focus on the case where $\lambda\gg \mu$, so that there is a separation between the timescale of
relaxation for the $x$ and $y$ variables. As such, $x$ is a `slow' variable, and is a natural
candidate for a reaction coordinate of the system; as in Section~\ref{sec:prob_formulation},
we therefore consider
\begin{equation*}
  \mbf{F}(\mbf{x}) = \Phi\mbf{x}\qquad\text{where}\qquad\Phi :=
  \left(\begin{array}{cc}
    1 & 0
  \end{array}\right).
\end{equation*}

\begin{figure}
\centering
\subfigure[$\mu =2, \lambda =20, \omega= 10, \tau =2$]{
\includegraphics[ width=0.45 \textwidth,height=4.8 cm]{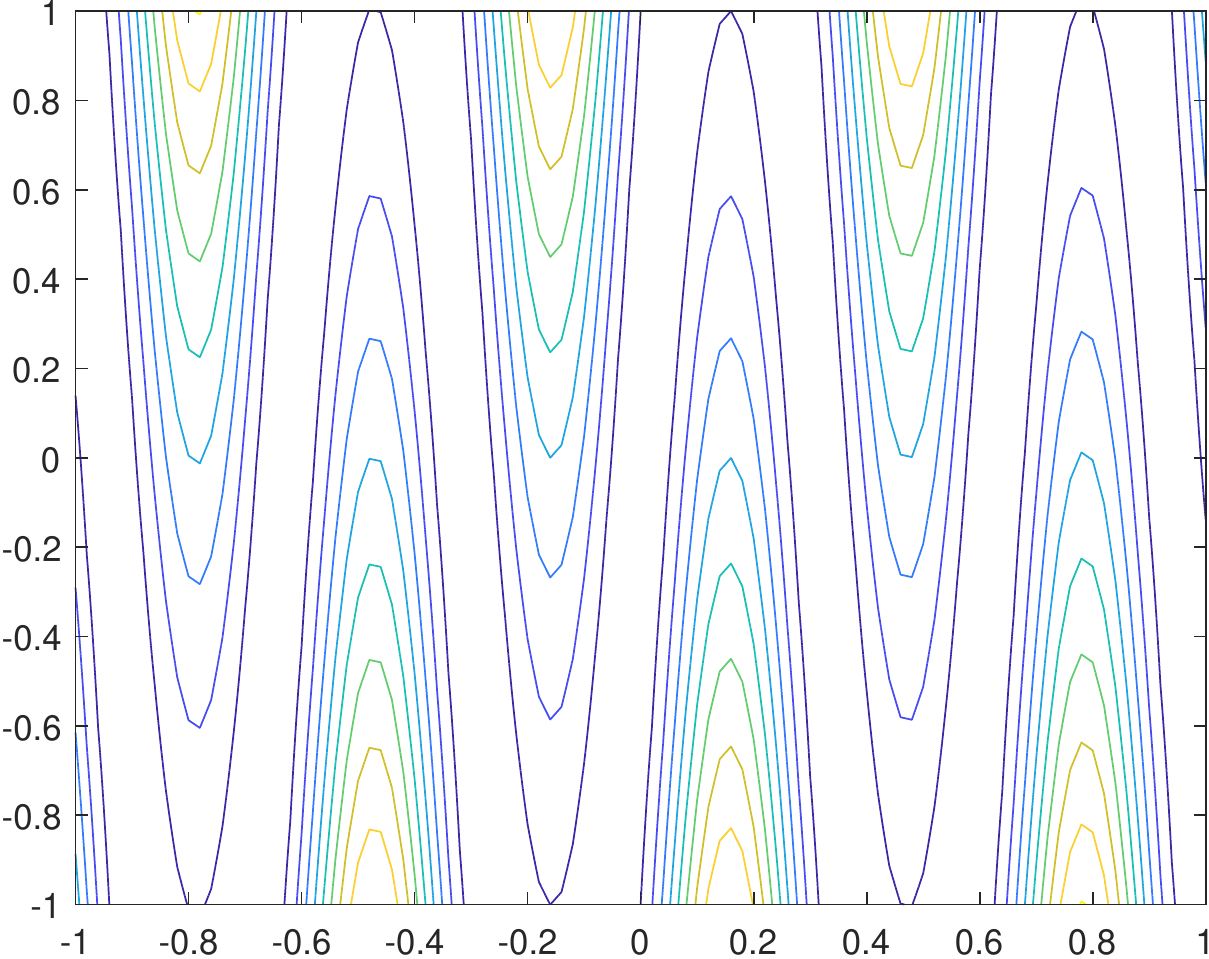}}
\;
\subfigure[$\mu =2, \lambda =20, \omega= 4, \tau =0.2$ ]{
\includegraphics[ width=0.45\textwidth,height=4.8 cm]{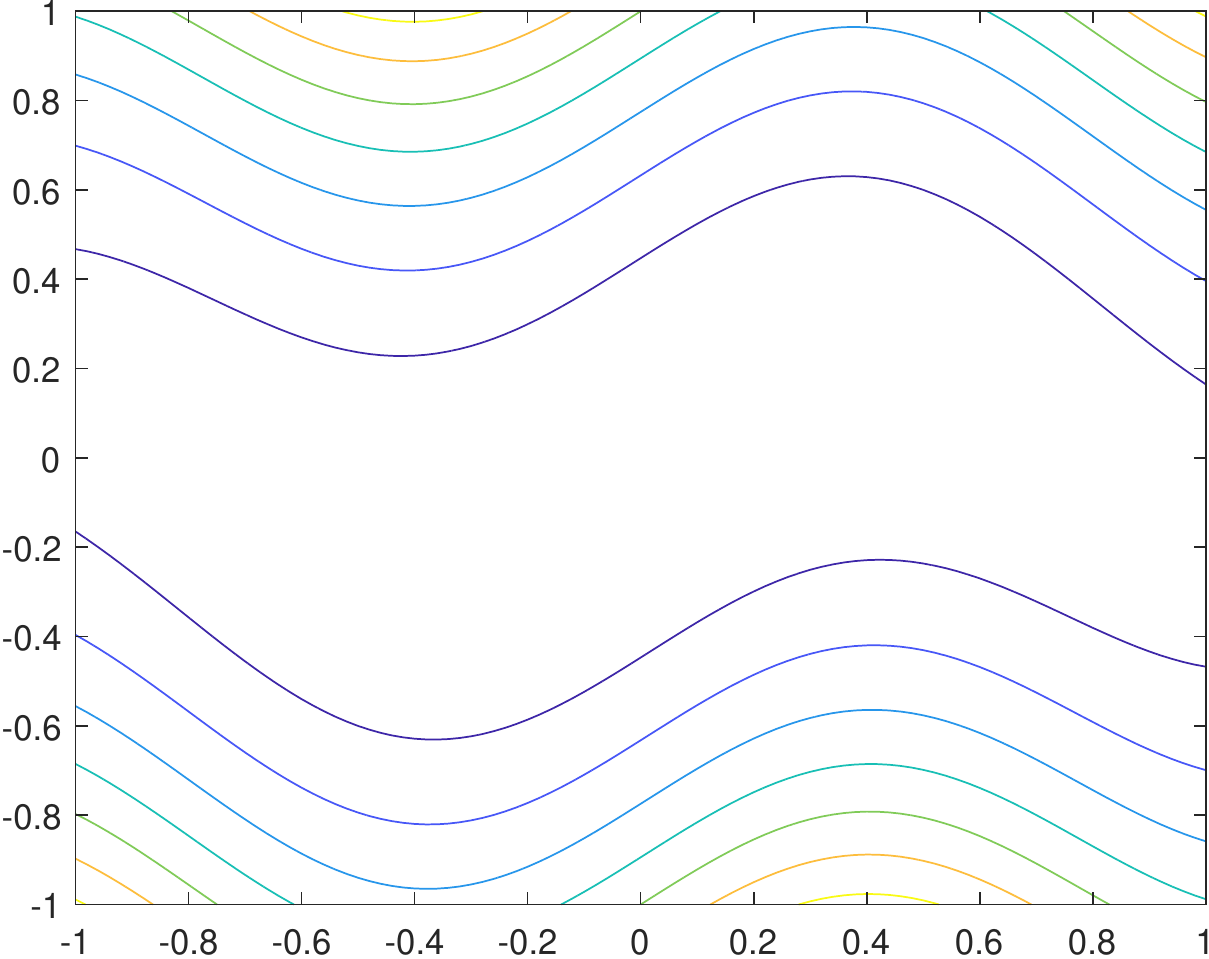}}
\vspace{- 10 pt}
\caption{Contour plots of the potential energy $V(x,y)$ defined in \eqref{potential_eq}: $\lambda/\mu \gg 1$
  and $\omega \tau$ is chosen to be $\gg 1$ for figure (a) and  $< 1$ for figure (b), respectively.
  }\label{Fig:energy_contour}
\end{figure}

The second term in \eqref{potential_eq}, i.e. $\frac{\lambda}{2}\big(\tau \sin(\omega x)-y\big)^2$,
has been chosen to emulate a form of free energy barrier to the dynamics, since when $\tau\sim 1$,
$\omega\gg1$ and $\beta\gg1$, we expect trajectories of the dynamics to remain close to the
manifold $y=\sin (\omega x)$; see Figure~\ref{Fig:energy_contour} for representations of different
energy landscapes.
As a consequence, we expect that as $\tau$ and $\omega$ increase with $\frac\lambda\mu\gg 1$, the
dynamics should to take progressively longer to approach a fixed neighbourhood of the global
equilibrium at $0$ from a generic initial condition, since the dynamics must effectively `travel
further' along the meandering valley in the potential to get there.

\subsection{Derivation of approximate dynamics}
\label{sec:Approx}
Under the assumptions described above, we compute $\calP\calL \mbf{F}$ and $\calQ\calL \mbf{F}$,
and use these to derive approximations of the terms involved in \eqref{MZ_formulation}.
\begin{enumerate}
\item{\emph{Computation of $\calP\calL \mbf{F}$.}
The effective potential $\calS(x)$ defined in  \eqref{effective_potential} is equal
to
\begin{equation}\label{eff_pot_ex}
\begin{split}
\calS(x)= & -\frac{1}{\beta}\log\left( \int_{\mbb{R}} e^{-\beta V(x,y)} dy \right)
=\frac{\mu}{2}x^2+\mathsf{const},
\end{split}
\end{equation}
as the unresolved variable $y$ follows normal distribution $y\sim N(\tau \sin (\omega x), \frac{1}{\beta \lambda})$, and hence
\begin{equation}\label{eff_dyn_eq}
\calP\calL \mbf{F}=
\left(
\begin{array}{c}
-\mu x\\
0
\end{array}
\right).
\end{equation}
}
\item\emph{Computation of $\calQ\calL \mbf{F}$.}
Clearly $\calQ\calL \mbf{F}=\calL \mbf{F}-\calP\calL \mbf{F}$, so using \eqref{eff_dyn_eq}, we have
\[
\calQ\calL \mbf{F}
=
\calQ\calL
\left(
\begin{array}{c}
x\\
y
\end{array}
\right)
=
\left(
\begin{array}{c}
-\lambda \tau \omega\left(\tau \sin (\omega x)-y\right)\,\cos(\omega x)\\
-\lambda\left(y-\tau \sin (\omega x)\right)
\end{array}
\right).
\]

\item\emph{Approximation of $\calM_s$.}
  Recalling the definition of the memory function $\calM_s(\mbf{F})$ from \eqref{memory_eq2},
  we must compute or otherwise approximate the expression $ \calQ\calL e^{s\calQ\calL} \mbf{F}
  \otimes  \calQ\calL \mbf{F}$. To do so, we define characteristic curves of the orthogonal
  dynamics,
  $\widetilde{\mbf{x}}_s=(\widetilde{x}_s,\widetilde{y}_s)=e^{s\calQ\calL}\widetilde{\mbf{x}}_0$, and find that
\begin{equation}\label{memory_dyn_eq1}
\begin{split}
& \calQ\calL e^{s\calQ\calL} \mbf{F}(\widetilde{\mbf{x}}_0)
\otimes  \calQ\calL \mbf{F}(\widetilde{\mbf{x}}_0)
=\Phi\dot{\widetilde{\mbf{x}}}_s
\otimes \Phi\dot{\widetilde{\mbf{x}}}_0\\
&\qquad= \lambda^2\tau^2\omega^2 \left(\tau\sin (\omega\widetilde{x}_s)-\widetilde{y}_s\right)\cos(\omega \widetilde{x}_s)
\cdot\,
\left(\tau\sin (\omega\widetilde{x}_0)-\widetilde{y}_0\right)\cos(\omega\widetilde{x}_0).
\end{split}
\end{equation}
We now change variable with the intention of linearizing, setting
$\widetilde{u}_s:=\widetilde{x}_s-\widetilde{x}_0$ and $\widetilde{v}_s:=\widetilde{y}_s-\tau\sin(\omega \widetilde{x}_0)$. Expressed in these new variables, the action of the orthogonal dynamics is equivalent to solving the ODE system
\begin{align*}
\left(
\begin{array}{c}
\dot{\widetilde{u}}_s\\
\dot{\widetilde{v}}_s
\end{array}
\right)
&=
\left(
\begin{array}{c}
-\lambda \tau \omega\left(\tau \sin (\omega(\widetilde{x}_0+\widetilde{u}_s))-\tau\sin(\omega \widetilde{x}_0)-\widetilde{v}_s\right)\,\cos\big(\omega (\widetilde{x}_0+\widetilde{u}_s)\big)\\
-\lambda\left(\tau\sin(\omega \widetilde{x}_0)+\widetilde{v}_s-\tau \sin \big(\omega (\widetilde{x}_0+\widetilde{u}_s)\big)\right)
\end{array}\right).
\end{align*}
Since we have assumed that $\lambda\gg \mu$, given knowledge of $\widetilde{x}_0$ alone, we expect that initial conditions for the orthogonal dynamics to be concentrated near
$(\widetilde{x}_0,\widetilde{y}_0) = (\widetilde{x}_0,\tau\sin(\omega \widetilde{x}_0))$, and so
linearising on this basis, we obtain
\begin{equation}\label{eq:linearisedOD}
\left(
  \begin{array}{c}
\dot{\widetilde{u}}_s\\
\dot{\widetilde{v}}_s
\end{array}
\right)
  =\lambda\left(
\begin{array}{cc}
  -\tau^2\omega^2\cos^2(\omega \widetilde{x}_0) & \tau \omega \cos(\omega \widetilde{x}_0) \\
  \tau\omega \cos(\omega \widetilde{x}_0) & -1
\end{array}\right)
\left(\begin{array}{c} \widetilde{u}_s\\ \widetilde{v}_s \end{array}\right)+\mathcal{O}(\widetilde{u}_s^2,\widetilde{v}_s^2,\widetilde{u}_s\widetilde{v}_s).
\end{equation}
Noting in particular that $\widetilde{u}_0=0$ since $\widetilde{u}_s=\widetilde{x}_s-\widetilde{x}_0$, and neglecting higher--order terms in \eqref{eq:linearisedOD}, the solution is approximately
\begin{align*}
  \left(
  \begin{array}{c}
\widetilde{u}_s\\
\widetilde{v}_s
\end{array}
\right)
  &\approx \frac{1}{1+\tau^2\omega^2\cos^2(\omega \widetilde{x}_0)}\left(\begin{array}{c}
   \widetilde{v}_0\tau\omega\cos(\omega \widetilde{x}_0) \\
   \widetilde{v}_0\tau^2\omega^2\cos^2(\omega\widetilde{x}_0)
  \end{array}\right)\\
  &\qquad\qquad+
\frac{e^{-\lambda\big(1+\tau^2\omega^2\cos^2(\omega\widetilde{x}_0)\big) s}}{1+\tau^2\omega^2\cos^2(\omega\widetilde{x}_0)}\left(\begin{array}{c}
   -\widetilde{v}_0\tau\omega\cos(\omega \widetilde{x}_0) \\
   \widetilde{v}_0
\end{array}\right),
\end{align*}
and therefore
\begin{equation}\label{eq:ApproxODVelocity}
  \dot{\widetilde{x}}_s = \dot{\widetilde{u}}_s \approx \lambda\tau\omega \widetilde{v}_0\cos(\omega \widetilde{x}_0) e^{-\lambda\big(1+\tau^2\omega^2\cos^2(\omega \widetilde{x}_0)\big) s}.
\end{equation}

When conditioning on knowledge of $\widetilde{x}_0$, it follows that
$\widetilde{y}_0\sim\mathcal{N}(\tau\sin(\omega \widetilde{x}_0),\frac{1}{\lambda\beta})$, and therefore $\widetilde{v}_0\sim\mathcal{N}(0,\frac{1}{\lambda\beta})$: the memory kernel is therefore approximately
\begin{equation}\label{approx_memory}
  \begin{aligned}
    \calM_s(\mbf{h}) &= \beta\,\mbb{E}[\calQ\calL e^{s\calQ\calL}\mbf{F}\otimes\calQ\calL\mbf{F}\,|\,\mbf{F}(\mbf{x})=\mbf{h}]\\
    &\approx \beta\int_{\mbb{R}}\dot{\widetilde{u}}_s(\mbf{h},\widetilde{y}_0)\otimes
    \dot{\widetilde{u}}_0(\mbf{h},\widetilde{y}_0)
    e^{-\beta\lambda\left(\widetilde{y}_0-\tau\sin(\omega\mbf{h})\right)^2}d\widetilde{y}_0\\
    &= \lambda\tau^2\omega^2 \cos^2(\omega \mbf{h})e^{-\lambda\big(1+\tau^2\omega^2\cos^2(\omega \mbf{h})\big) s},
  \end{aligned}
\end{equation}
where in the above formula, $\dot{\widetilde{u}}_0=\lim_{s\to 0+}\dot{\widetilde{u}}_s$.
Recalling the form of \eqref{MZ_formulation}, we note that we must also approximate the divergence
of $\calM_s$; using the expression derived in \eqref{approx_memory}, in this case we obtain
\begin{equation}\label{approx_memory_div}
  \mathsf{div}\calM_s(\mbf{h}) \approx -\lambda\tau^2\omega^3\sin(2\omega\mbf{h})
  \left(1-\lambda\tau^2\omega^2\cos^2(\omega\mbf{h})s\right)
  e^{-\lambda\big(1+\tau^2\omega^2\cos^2(\omega \mbf{h})\big) s}.
\end{equation}

\item\emph{Approximation of memory integral.}
  Our next step is to approximate the first integral term involving the memory kernel in
  \eqref{MZ_formulation}.
  Noting the form of \eqref{approx_memory}, we see that this is an integral of exponential type,
  and therefore we apply the method of steepest descents to derive an approximation.

  We note that the exponential term in the approximate expressions \eqref{approx_memory} and
  \eqref{approx_memory_div} are maximal when $s=0$; if $\mbf{h}_{t-s}\approx \mbf{h}_t$ for small
  $s$, we obtain
\begin{equation*}
  \begin{aligned}
    \int_{0}^{t} \calM_s(\mbf{h}_{t-s})\cdot \nabla\calS(\mbf{h}_{t-s}) ds
    &\approx \int_0^t \lambda\tau^2\omega^2 \cos^2(\omega \mbf{h}_{t-s})\,\mu \mbf{h}_{t-s}
    e^{-\lambda(1+\tau^2\omega^2\cos^2(\omega \mbf{h}_{t-s})) s} ds\\
    &\approx\int_0^\infty \lambda\tau^2\omega^2 \cos^2(\omega \mbf{h}_{t})\,\mu \mbf{h}_{t}e^{-\lambda(1+\tau^2\omega^2\cos^2(\omega \mbf{h}_t)) s} ds\\
    &\approx \frac{\tau^2\omega^2 \cos^2(\omega \mbf{h}_{t})}{1+\tau^2\omega^2\cos^2(\omega \mbf{h}_t)}\mu \mbf{h}_{t},
\end{aligned}
\end{equation*}
and via a similar approximation for the divergence term, we have
\begin{equation*}
  \begin{aligned}
    &\int_{0}^{t}-\frac{1}{\beta}\mathsf{div}\left(\calM_s(\mbf{h}_{t-s})\right) ds\\
    &\qquad\approx \int_0^t \frac{1}{\beta}\lambda\tau^2\omega^3\sin(2\omega\mbf{h}_{t-s})
  \left(1-\lambda\tau^2\omega^2\cos^2(\omega\mbf{h}_{t-s})s\right)
  e^{-\lambda(1+\tau^2\omega^2\cos^2(\omega \mbf{h}_{t-s})) s} ds\\
  &\qquad \approx \int_0^\infty \frac{1}{\beta}\lambda\tau^2\omega^3\sin(2\omega\mbf{h}_{t})
  \left(1-\lambda\tau^2\omega^2\cos^2(\omega\mbf{h}_{t})s\right)
  e^{-\lambda(1+\tau^2\omega^2\cos^2(\omega \mbf{h}_{t})) s} ds\\
  &\qquad\approx \frac{1}{\beta}\frac{\tau^2\omega^3\sin(2\omega\mbf{h}_{t})}
  {(1+\tau^2\omega^2\cos^2(\omega \mbf{h}_{t}))^2}.
\end{aligned}
\end{equation*}
Combining these approximations, we obtain an approximation of the memory contributions in \eqref{MZ_formulation} as
\begin{equation}\label{eq:approx_memory_integral}
\begin{split}
&\int_{0}^{t} \calM_s(\mbf{h}_{t-s})\cdot \nabla \calS(\mbf{h}_{t-s})-\frac{1}{\beta}\mathsf{div}\,\calM_s (\mbf{h}_{t-s})ds\\
&\qquad\qquad\approx \frac{\tau^2\omega^2 \cos^2(\omega \mbf{h}_{t})}{1+\tau^2\omega^2\cos^2(\omega \mbf{h}_t)}\mu \mbf{h}_{t}+\frac{1}{\beta}\frac{\tau^2\omega^3\sin(2\omega\mbf{h}_{t})}
  {(1+\tau^2\omega^2\cos^2(\omega \mbf{h}_{t}))^2}.
\end{split}
\end{equation}
Notably, if $\tau\sim 1$, the second term is negligible compared with the first both when
$\omega\ll 1$, and when $\omega\gg 1$. We will therefore discard the latter term in these cases.
\item \emph{Fluctuating force.} Above, we have shown that the memory
  kernel can be approximated as
  \begin{equation*}
    \calM_s(\mbf{h}) \approx \frac{\tau^2\omega^2\cos^2(\omega\mbf{h})}{1+\tau^2\omega^2\cos^2(\omega\mbf{h})}\delta_0(s).
  \end{equation*}
  Since $\calM_s$ is the autocovariance of $\calF_t$ multiplied by $\beta$,
  and the autocovariance of the white noise already present is
  $\beta^{-1}\delta_0(s)$, it is natural, in view of the Fluctuation--Dissipation Theorem,
  to define a new stochastic forcing which has an autocovariance that is the sum of these two
  contributions, i.e.
  \begin{equation*}
    d\calF_t+\sqrt{2\beta^{-1}}d\mbf{B}_t\approx\sqrt{\frac{2\beta^{-1}}{1+\tau^2\omega^2\cos^2(\omega\mbf{h})}} d\mbf{B}_t
  \end{equation*}
\end{enumerate}
Combining the approximations above, in the case where $\omega\gg1$, we obtain the closed--form approximate equation
\begin{equation}\label{ApproxDyn}
  d\mbf{h}_t=-\frac{\mu\mbf{h}_t}{1+\tau^2\omega^2\cos^2(\omega \mbf{h}_t)}dt
  +\sqrt{\frac{2\beta^{-1}}{1+\tau^2\omega^2\cos^2(\omega\mbf{h}_t)}} d\mbf{B}_t
\end{equation}
Notably, the drift term in this equation is independent of $\beta$.

\subsection{Other choices of approximate dynamics}
The derivation of the effective dynamics \eqref{eq:approx_memory_integral} was informed by
the Mori--Zwanzig formalism, but other choices could be made, and may be more appropriate in
other circumstances.

\begin{enumerate}
\item \emph{Discarding memory and fluctuating force.}
  In \cite{Legoll2010a}, the authors consider another choice of effective dynamics, which in our
  setting, amounts to considering
  \begin{equation}\label{FredTony}
    d\boldsymbol{\xi}_t = -\Phi\Phi^T\nabla\calS(\boldsymbol{\xi}_t)dt+\sqrt{2\beta^{-1}\Phi\Phi^T}
    d\mbf{B}_t.
  \end{equation}
  This is equivalent to \eqref{MZStochastic} where the memory and fluctuating force terms have been
  neglected entirely. For the evolution of the mean $\mbf{h}_t$, this choice of dynamics yields
  \begin{equation}\label{NoMemory}
    \partial_t\mbf{h}_t = -\Phi\Phi^T \nabla \calS(\mbf{h}_t)dt
  \end{equation}
  as the effective equation for the observable we consider.
  The authors have proved error bounds on the time marginals of the resulting probability
  distribution when compared the true dynamics captured by \eqref{Liouville_eq2}; applying
  \cite[Proposition~3.1]{Legoll2010a} to our case gives the bound
  \begin{equation}\label{RelEntBound}
    H(\psi_t|\phi_t)\leq \frac{\beta\tau^2\omega^2}{4} \left[H\left(\psi_0\middle|\mu\right)-H\left(\psi_t\middle|\mu\right)\right],
  \end{equation}
  where:
  \begin{enumerate}
  \item $H(\mu|\nu)$ is the relative entropy of a measure $\mu$ with respect to $\nu$, i.e.
    \begin{equation*}
      H(\mu|\nu):= \int\log\left(\frac{d\mu}{d\nu}\right)d\mu;
    \end{equation*}
  \item $\mu$ is the Gibbs measure;
  \item $\psi_t$ is the distribution of the `true' dynamics at time $t$; and
  \item $\phi_t$ is the distribution of solutions to \eqref{FredTony} at time $t$.
  \end{enumerate}
  Clearly, the constant in \eqref{RelEntBound} is large when $\tau\omega \gg 1$;
  this reflects the fact that neglecting the memory in this case is not sufficient to accurately
  capture the dynamical properties of the system, and a more sophisticated approach is needed.
  \medskip
\item \emph{A more na\"ive memory approximation}.
  To highlight the need to conduct dynamical sampling to approximate $\calM_s$ correctly,
  we remark that the approximation of the memory terms obtained in
  \eqref{eq:approx_memory_integral} is notably \emph{not} the same as simply choosing to
  approximate
  \begin{align*}
    \calM_s(\mbf{h})
    &\approx\widetilde{\calM}_s(\mbf{h}):=
      \beta\,\mbb{E}[\calQ\calL\mbf{x}\otimes\calQ\calL\mbf{x}\,|\,\mbf{F}(\mbf{x})
      =\mbf{h}]\,\delta_{0}(s),\\
    &=\left(\sqrt{\frac{\lambda\beta}{2\pi}}\int \beta\,\lambda^2\tau^2\omega^2(\tau
      \sin(\omega \mbf{h})-y)^2\cos^2(\omega \mbf{h})
      e^{-\frac12\beta\lambda(\tau \sin(\omega \mbf{h})-y)^2}dy\right)\,\delta_0(s)\\
    &=\lambda\tau^2\omega^2\cos^2(\omega \mbf{h})\delta_0(s).
  \end{align*}
  Using $\widetilde{\calM}_s$ would result in the approximate dynamics
  \begin{equation}\label{memory_naive}
    \begin{split}
      \partial_t\mbf{h}_t&=\big(\lambda\tau^2\omega^2 \cos^2(\omega \mbf{h}_{t})-1\big)\mu \mbf{h}_t
      +\frac{\lambda}{\beta}\tau^2\omega^3\sin(2\omega\mbf{h}_{t}).
    \end{split}
  \end{equation}
  Since we have chosen $\lambda\gg 1$, we see that this approximation will yield qualitatively
  different dynamics to both the true dynamics for $\mbf{F}$, \eqref{Liouville_eq2}, and the
  approximate dynamics given by \eqref{eq:approx_memory_integral}; we investigate this numerically
  in Section~\ref{sec:numerics}.
\end{enumerate}
The remarks above suggests that in general, careful dynamical sampling of the memory kernel is
required to accurately capture the interaction between chosen reaction coordinates and the
neglected degrees of freedom.

\section{Numerical simulations}\label{sec:numerics}
In this section, we conduct a numerical study of the various choices of approximate effective
dynamics for the observable in the toy example considered in Section~\ref{sec:Approx}.
In particular, we inspect the validity of some of
the approximations made, and compare the different choices of effective dynamics described in
Section~\ref{sec:benchmark}.


\subsection{Investigation of the memory kernel}
The derivation of the effective dynamics for $\mbf{h}_t$ made in Section~\ref{sec:Approx} relies
crucially upon a series of approximations to the memory integral in \eqref{MZ_formulation}; we
therefore first consider the validity of these assumptions.

\begin{enumerate}
\item\emph{Sampling of the memory kernel.}
  To understand the error committed by using the approximation of the memory kernel given in
  \eqref{approx_memory}, we computed the memory kernel by statistically sampling trajectories of
  the orthogonal dynamics. The result of this simulation for 2 different cases when $\tau\omega$
  is of different sizes are shown in Figure~\ref{Fig:memory}. The empirical memory exhibits
  monotone decay in time for all values of $x_0$ considered, and decay is more rapid when
  $\tau\omega$ is large, in agreement with \eqref{approx_memory}.
\begin{figure}[htp!]
\centering
\subfigure[ $\omega =10 $, $\tau=2$ ]{
\includegraphics[ width=0.45\textwidth,height=5 cm]{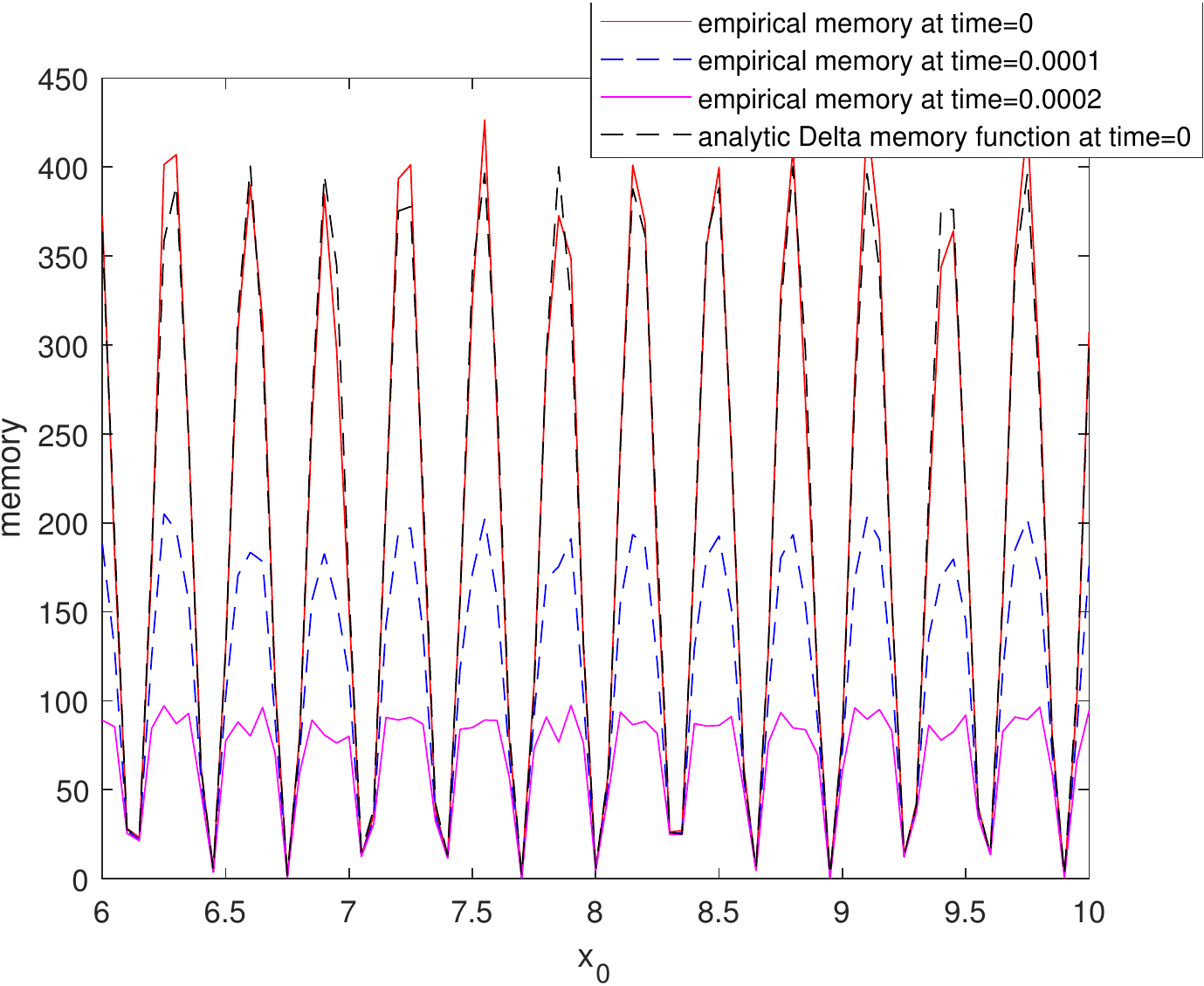}}
\subfigure[$\omega =4 $, $\tau=.2$]{
\includegraphics[ width=0.45\textwidth,height=5 cm]{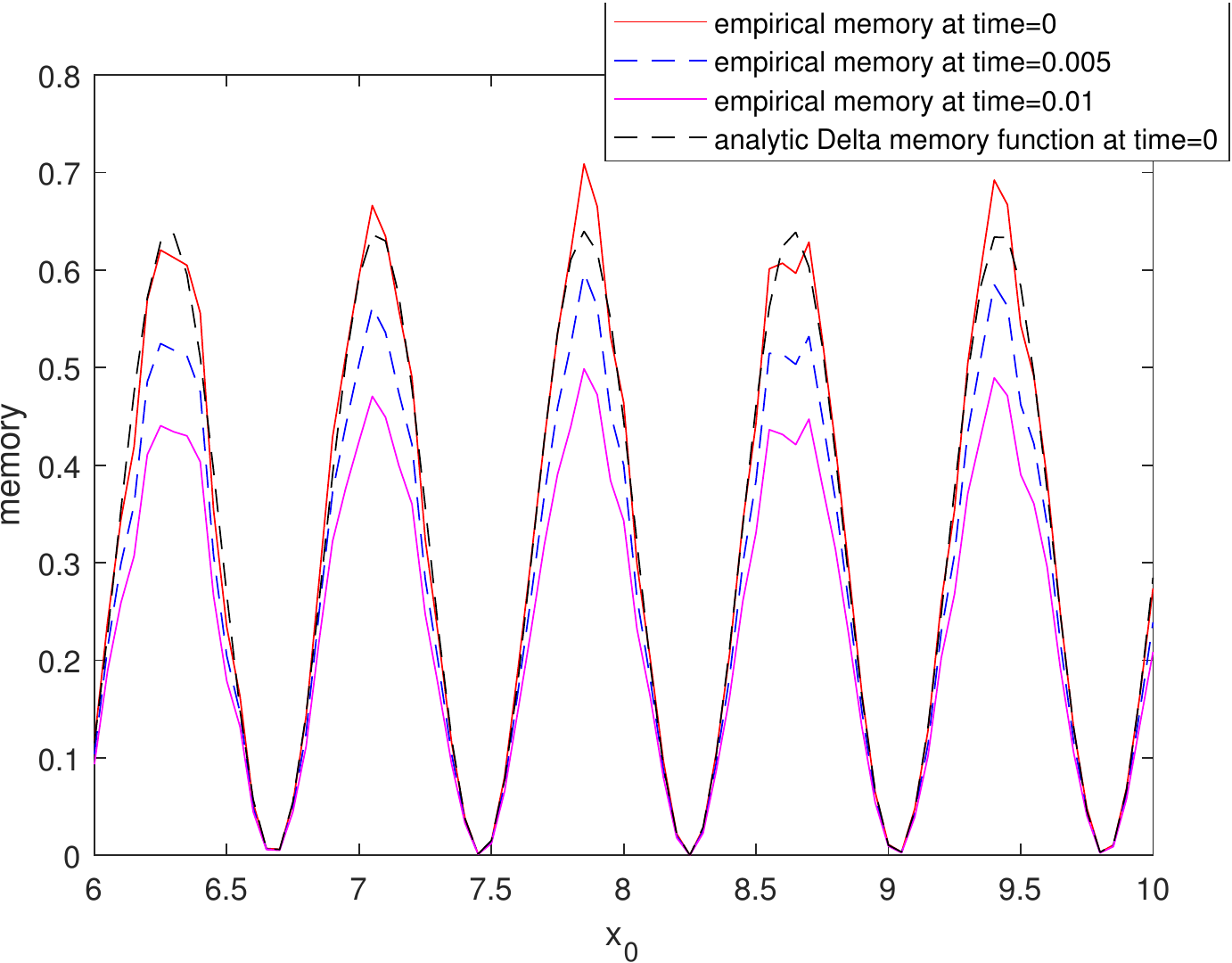}}
\vspace{- 10 pt}
\caption{Comparison of empirical memory kernel and approximate form given
  in \eqref{approx_memory}. In both cases $\lambda=20$ and $\mu=2$.\label{Fig:memory}}
\vspace*{\floatsep}
\centering
\subfigure[ $|\cos(\omega x_0)|= 1$ ]{
  \includegraphics[ width=0.45\textwidth,height=4.8 cm]{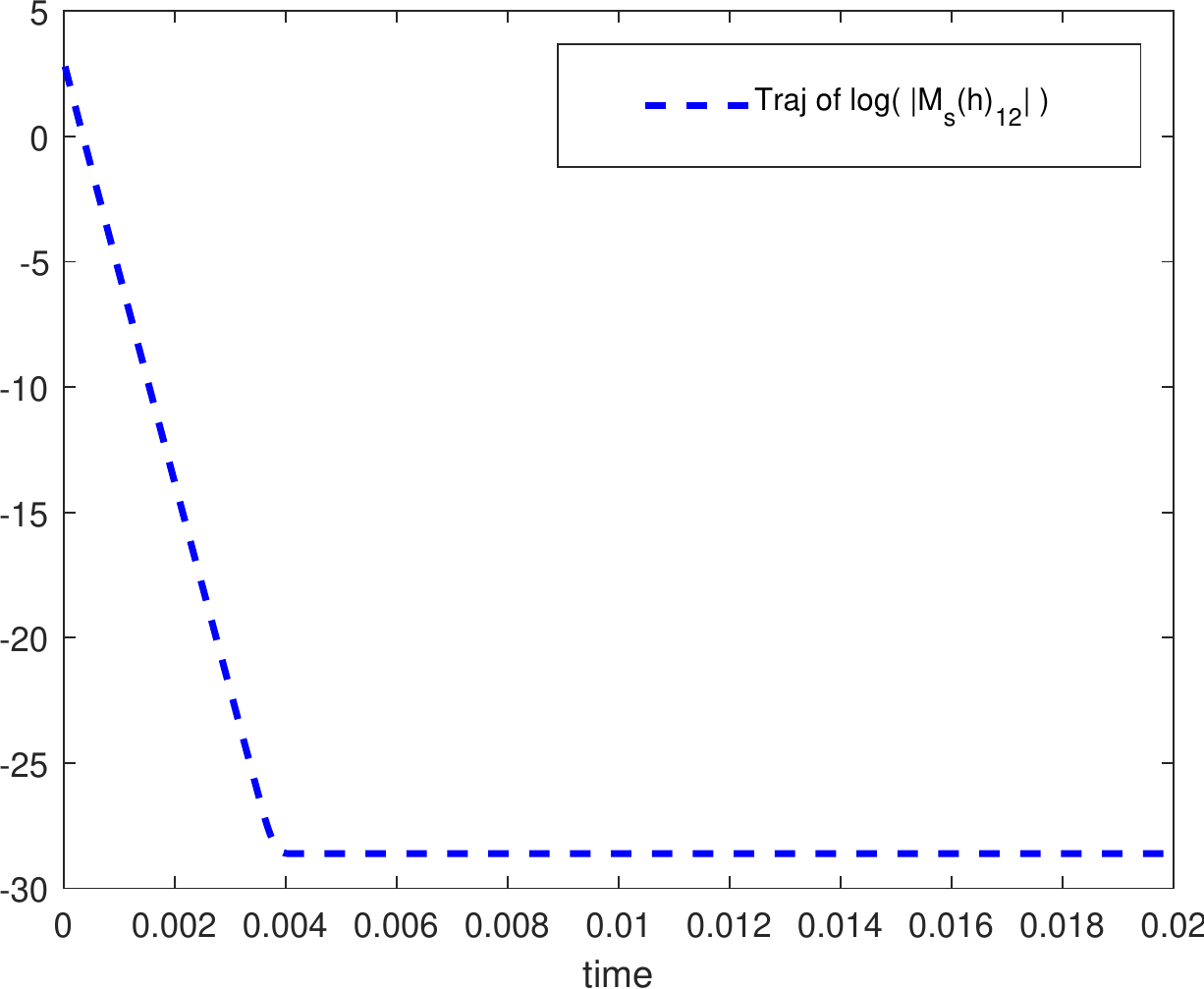}
}
\subfigure[$\cos(\omega x_0)= 0$]{
  \includegraphics[ width=0.45\textwidth,height=4.8 cm]{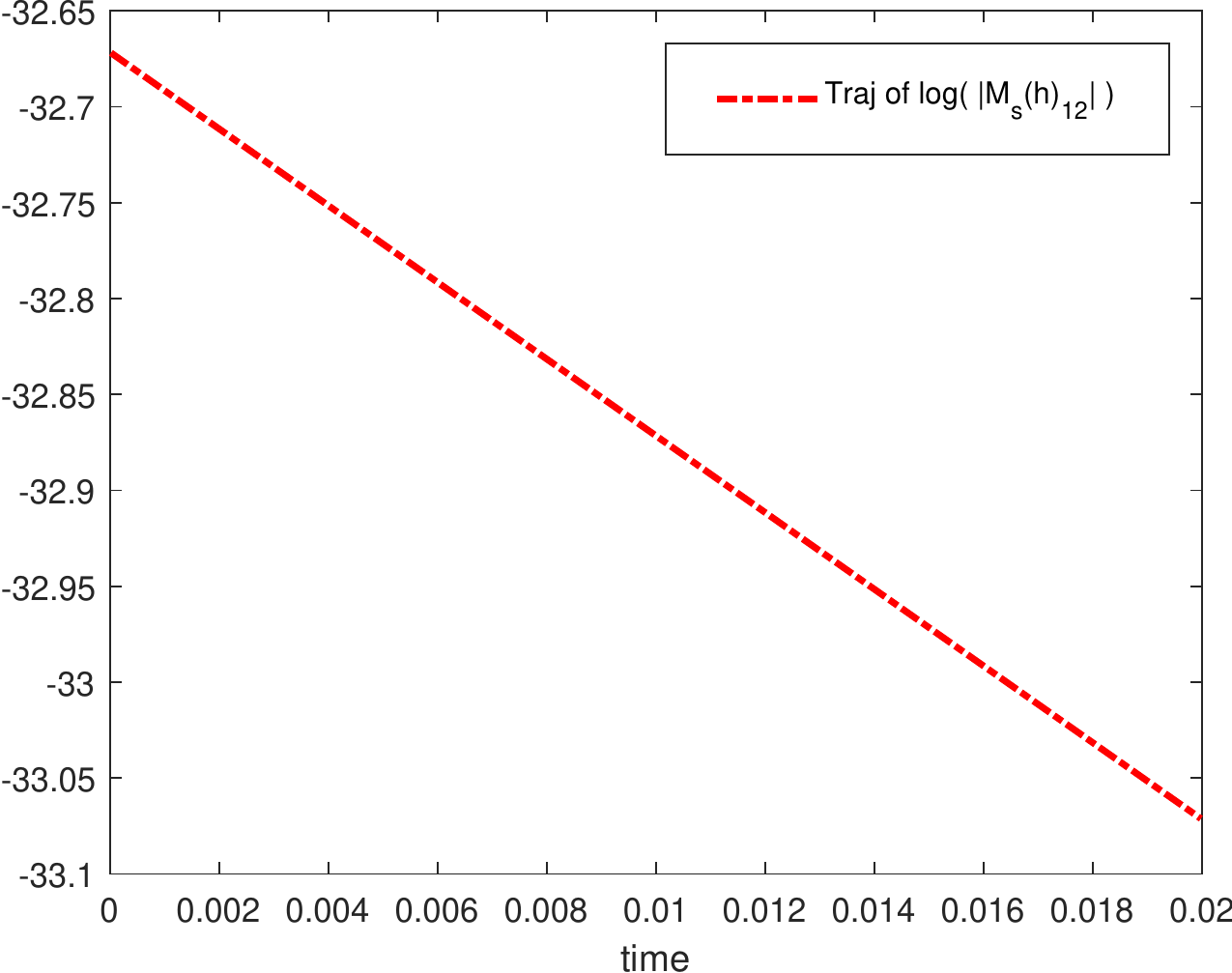}
}
\vspace{- 1 pt}
\caption{Average trajectories of $\log\left|\mbf{M}_s(\mbf{h})_{12}\right|$ under different initial condition of $x_0$ (see \eqref{memory_matrix}) based on a sample of 2000 initial conditions of $y_0$. Initial conditions of $y_0$ satisfy
  $y_0\sim\mathcal{N}(\tau \sin(\omega x_0), \frac{1}{\beta\lambda})$, and parameters of $V$ are
  $\lambda=20$, $\mu=2$, $\tau=2$ and $\omega=10$.\label{Fig:memory_matrix}}
\end{figure}

\item\emph{The memory kernel and choice of reaction coordinate}.
  In practice, it is impractical to compute and then integrate over long trajectories of the memory
  term, so a good coarse-graining selector should lead to a fast decay of the memory kernel
  \cite{Guttenberg2013a}.

  With this in mind, in general we set $\Sigma=\sqrt{\Phi\Phi^T}$, and consider
  \begin{equation}\label{memory_matrix}
  \begin{aligned}
    \mbf{M}_s(\mbf{h})&:=\beta\mbb{E}\big[  e^{s\calQ\calL}\calQ\calL \mbf{x} \otimes  \calQ\calL \mbf{x}\,\big|\, \Phi \mbf{x}=\mbf{h}\big]\\
    &= \beta\mbb{E}\left[\left(
        \begin{array}{cc}
          e^{s\calQ\calL}\calQ\calL\Sigma^{-1}\Phi\mbf{x} \otimes  \calQ\calL \Sigma^{-1}\Phi\mbf{x}
          &e^{s\calQ\calL}\calQ\calL\Sigma^{-1}\Phi\mbf{x} \otimes  \calQ\calL \Psi\mbf{x}\\
          e^{s\calQ\calL}\calQ\calL\Psi\mbf{x} \otimes  \calQ\calL \Sigma^{-1}\Phi\mbf{x}
          &e^{s\calQ\calL}\calQ\calL\Psi\mbf{x} \otimes  \calQ\calL \Psi\mbf{x}
        \end{array}                                                                                                   \right)\,\middle|\, \Phi \mbf{x}=\mbf{h}\right]\\
    &= \left(
        \begin{array}{cc}
          \mbf{M}_s(\mbf{h})_{11}
          &\mbf{M}_s(\mbf{h})_{12}\\
          \mbf{M}_s(\mbf{h})_{21}
          &\mbf{M}_s(\mbf{h})_{22}
        \end{array}                                                                                                   \right).
  \end{aligned}
\end{equation}
We see that $\mbf{M}_s(\mbf{h})_{11} = \Sigma^{-1}\calM_s(\mbf{h})\Sigma^{-1}$, recalling the
definition of $\calM_S$ from \eqref{memory_eq2}. We also note that
if dynamical sampling of the orthogonal
dynamics is used to approximate $\calM_s$ in practice, all of the information needed to compute
$\mbf{M}_s$ is available.

Intuitively, the off--diagonal blocks of $\mbf{M}_s$ describe the correlation
between the action of the fluctuating force on the reaction coordinates and orthogonal variables,
and in particular, $\mbf{M}_s(\mbf{h})_{12}$ describes the influence of the unknown variables
$\Psi\mbf{x}$ at the current time on the dynamics of $\Phi\mbf{x}$ at later times.
We can therefore test the strength of the `coupling' between the reaction coordinates
and the other degrees of freedom by considering the magnitude of $\mbf{M}_s(\mbf{h})_{12}$;
a log plot is shown in Figure~\ref{Fig:memory_matrix} demonstrating very rapid decay of the
corresponding entry.

We note that the errors introduced by the delta approximation of memory kernel are accumulated over
time, so, if the $L^{\infty}$ norms of entire trajectories of $\mbf{M}_s(\mbf{h})_{12}$ is small, (that is when $\cos(\omega x_0)$ is equal to zero), then the memory kernel decays faster and the delta approximation \eqref{eq:approx_memory_integral} is more appropriate. 

\end{enumerate}
\subsection{Comparison of different effective dynamics}
In this section, we compare the various choices of effective dynamics described in
Section~\ref{sec:benchmark} with the true dynamics.
In order to compare with the `true' dynamics, given by \eqref{MZ_formulation}, we perform
statistical sampling on the full dynamics, given by \eqref{Liouville_eq2}.

\begin{enumerate}
\item\emph{Simulations without thermostat}.
  We first simulated the evolution of the mean value of the observable $\mbf{F}(\mbf{X}_t)$
  for various choices of dynamics without thermostat. For convenience, we recall the relevant
  governing equations are:
  \begin{align}
    \partial_t\mbf{h}_t&=-\frac{\mu\mbf{h}_t}{1+\tau^2\omega^2\cos^2(\omega \mbf{h}_t)}
                         \tag{\ref{ApproxDyn}}\\
    \partial_t\mbf{h}_t&=-\mu\mbf{h}_t\tag{\ref{NoMemory}}\\
    \partial_t\mbf{h}_t&=\big(\lambda\tau^2\omega^2 \cos^2(\omega \mbf{h}_{t})-1\big)\mu \mbf{h}_t
                         +\frac{\lambda}{\beta}\tau^2\omega^3\sin(2\omega\mbf{h}_{t}).
                         \tag{\ref{memory_naive}}
  \end{align}
  We fixed an initial condition $\mbf{h}_0=x_0$, and compared with the `true' dynamics without
  thermostat, i.e. we considered $\mbf{h}_t:=\mbb{E}[\mbf{F}(\mbf{X}_t)|\mbf{F}(\mbf{X}_0)=x_0]$
  by sampling
  \begin{equation}\label{NoThermostat}
      d\mbf{F}(\mbf{X}_t) = -\Phi \nabla V(\mbf{X}_t)  dt,\quad
      \text{where}\quad X_0=x_0\quad\text{and}\quad
      Y_0\sim \mathcal{N}\left(\tau \sin(\omega x_0),\tfrac{1}{\beta\lambda}\right).
  \end{equation}
  Figure~\ref{Fig:full_vs_MZ} shows the results of these simulations, showing that the trajectories
  of \eqref{ApproxDyn} and averages of \eqref{NoThermostat} match with each other very well,
  whereas the effective system \eqref{NoMemory} neglecting the memory contributions relaxes much faster, and the
  na\"ive choice of memory made to arrive at \eqref{memory_naive} yields incorrect behaviour.
  \medskip

\begin{figure}[htp!]
\centering
\subfigure[ $\tau=2$, $\omega=10$]{
\includegraphics[ width=0.4\textwidth, height=5 cm]{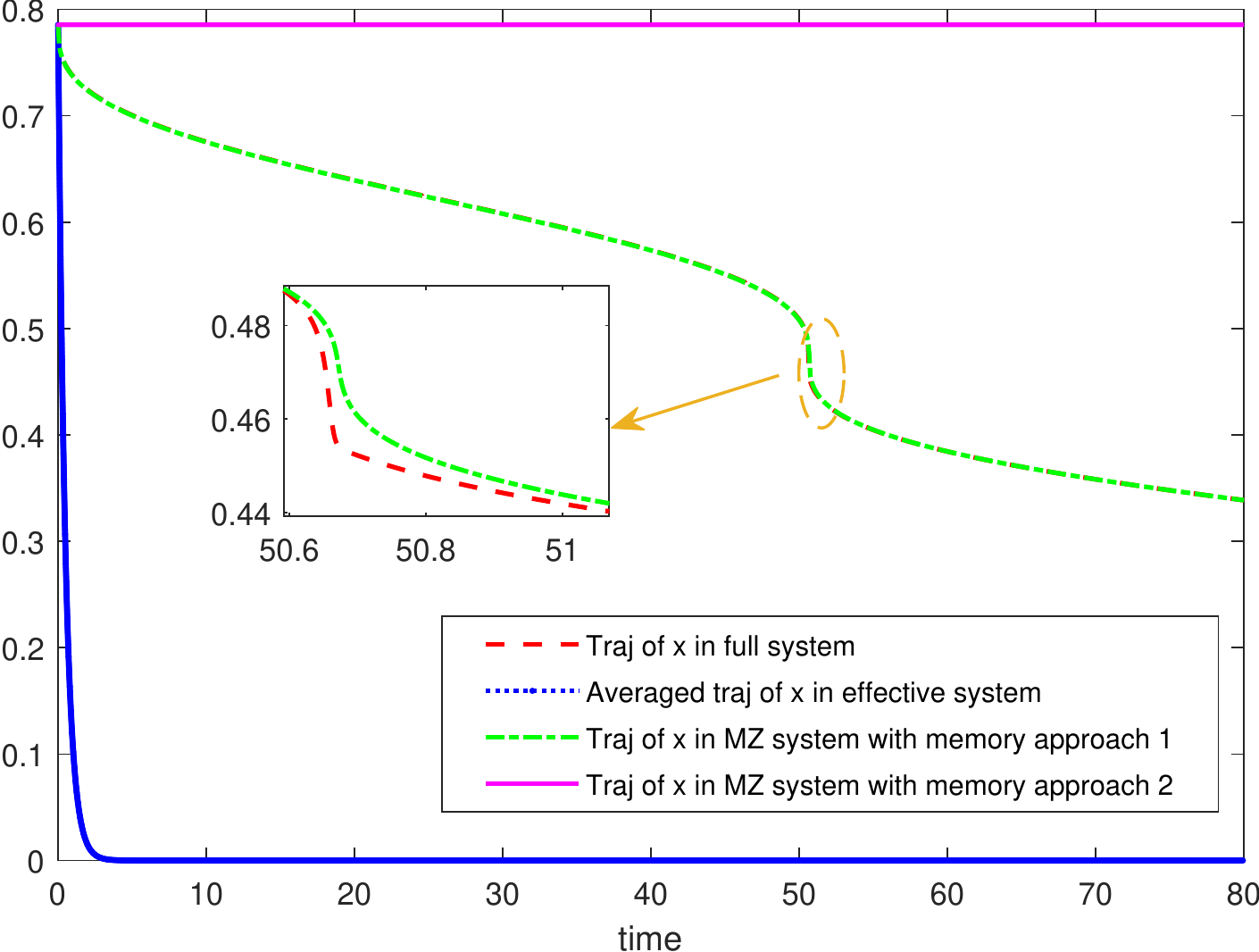}}
\quad
\subfigure[$\tau =0.2$, $\omega=4$]{
\includegraphics[ width=0.4\textwidth, height=5 cm]{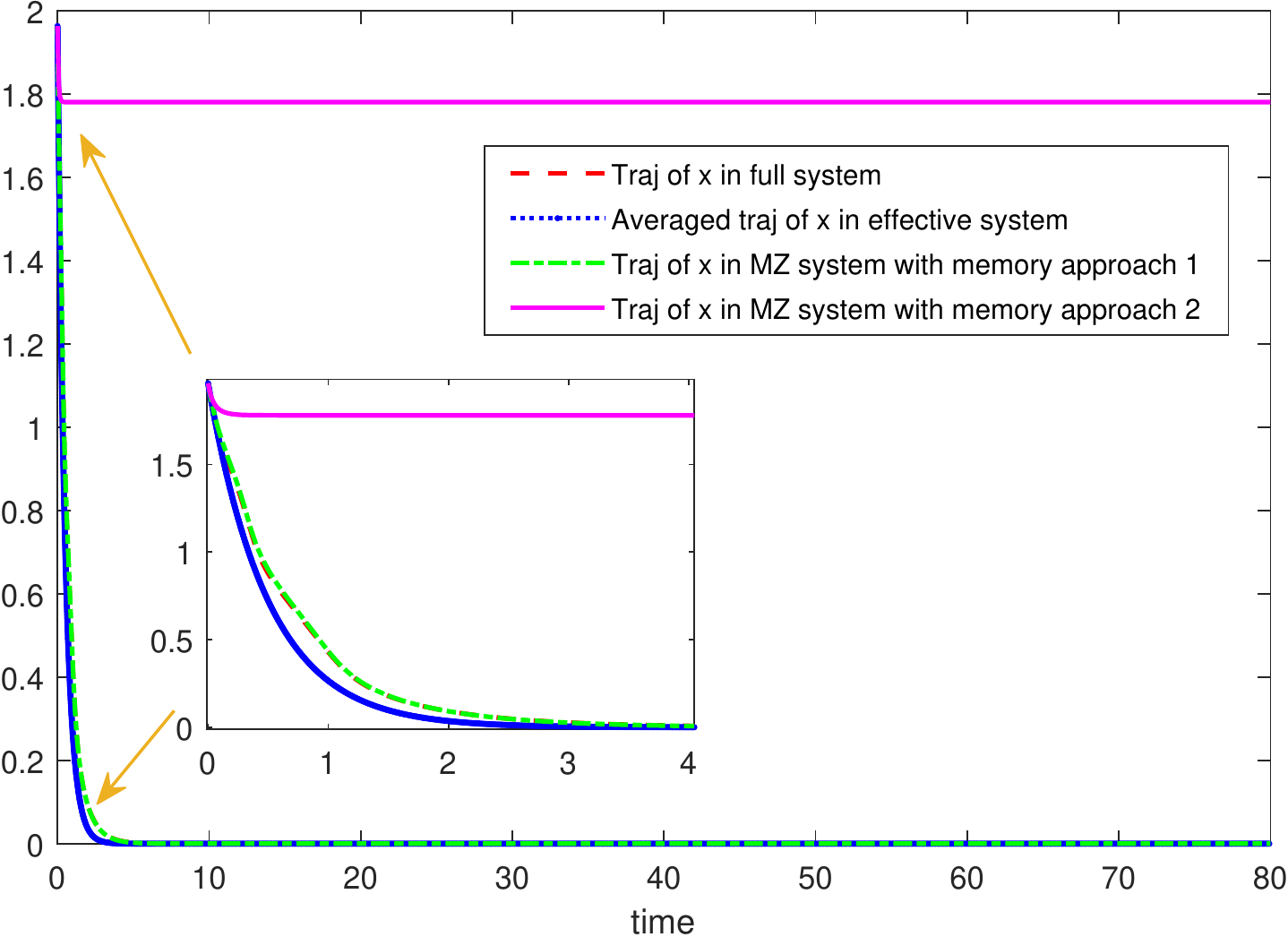}}
\vspace{- 0.2 pt}
\caption{Trajectories of $\mbf{h}_t$ evolving under \eqref{NoThermostat}, \eqref{NoMemory}, \eqref{ApproxDyn} (Approach 1) and \eqref{memory_naive} (Approach 2), respectively shown in red, green, blue and magenta.
  Parameters of $V$ are $\lambda=20$, $\mu=2$, and the time step was $\Delta t=10^{-5}$ and $T=80$.\label{Fig:full_vs_MZ}}
\vspace*{\floatsep}
\centering
\subfigure[ $\beta=1$]{
\includegraphics[ width=0.4\textwidth, height=5 cm]{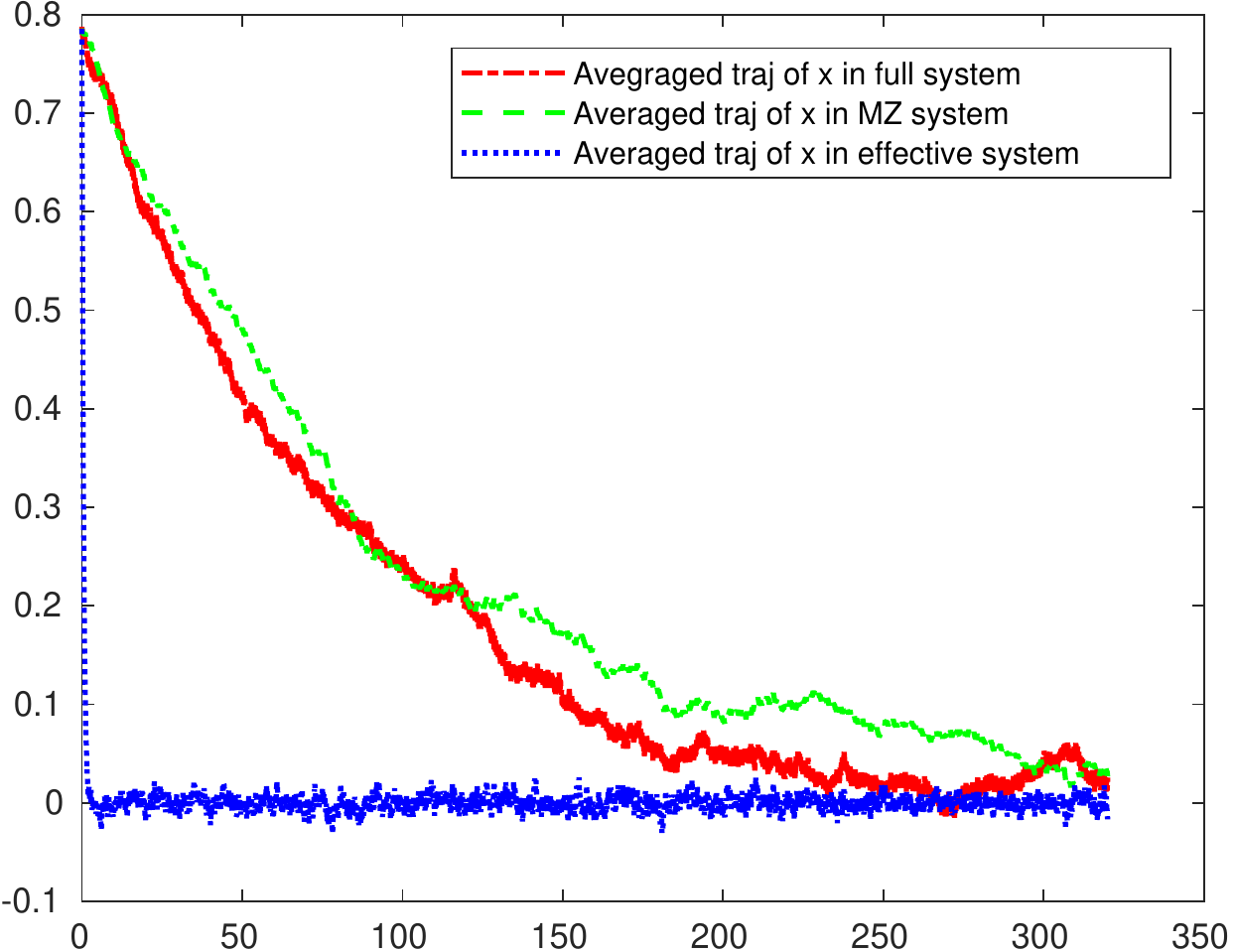}}
\subfigure[$\beta=10$]{
\includegraphics[ width=0.4 \textwidth, height=5 cm]{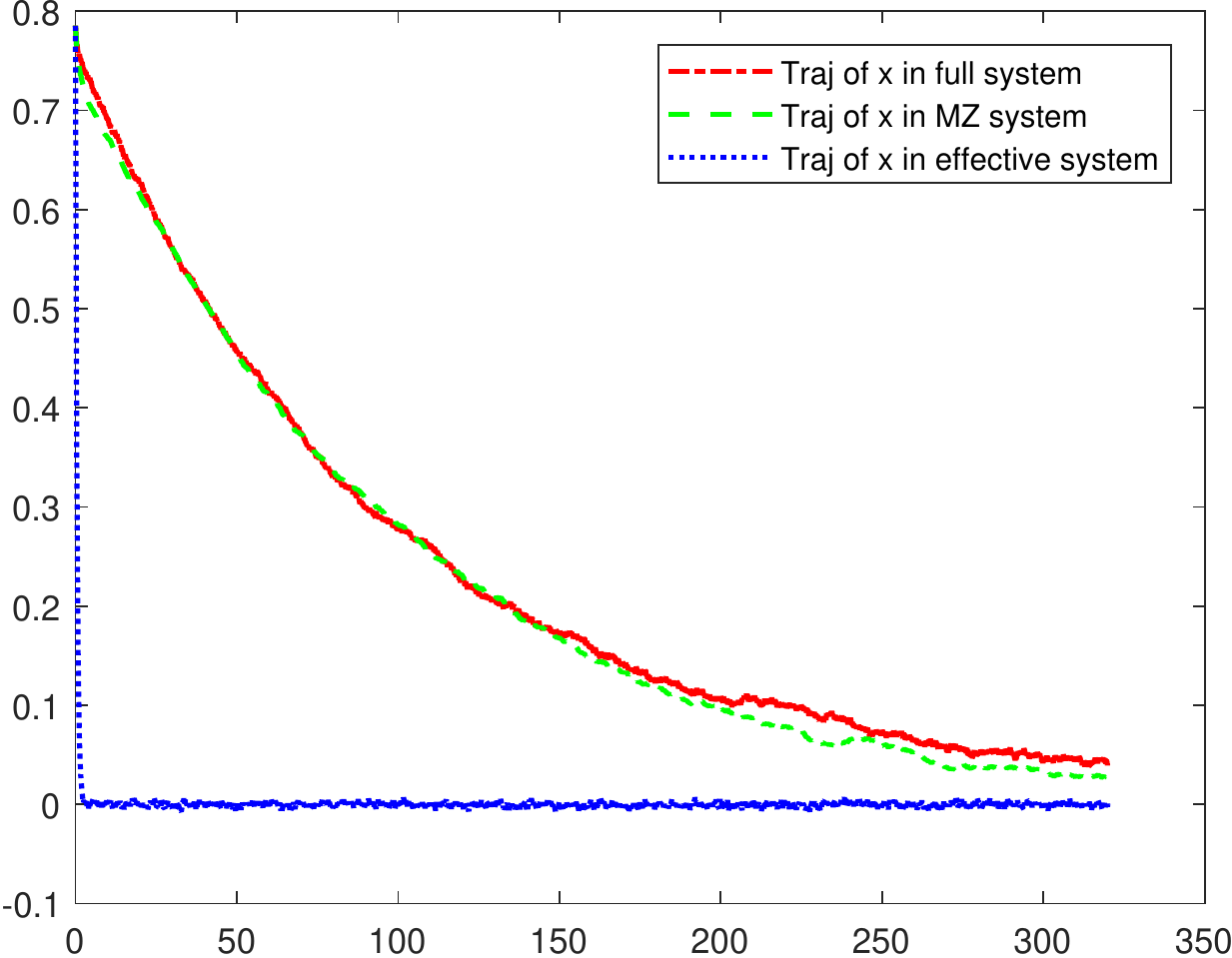}}
\subfigure[ $\beta=100$]{
\includegraphics[ width=0.4\textwidth, height=5 cm]{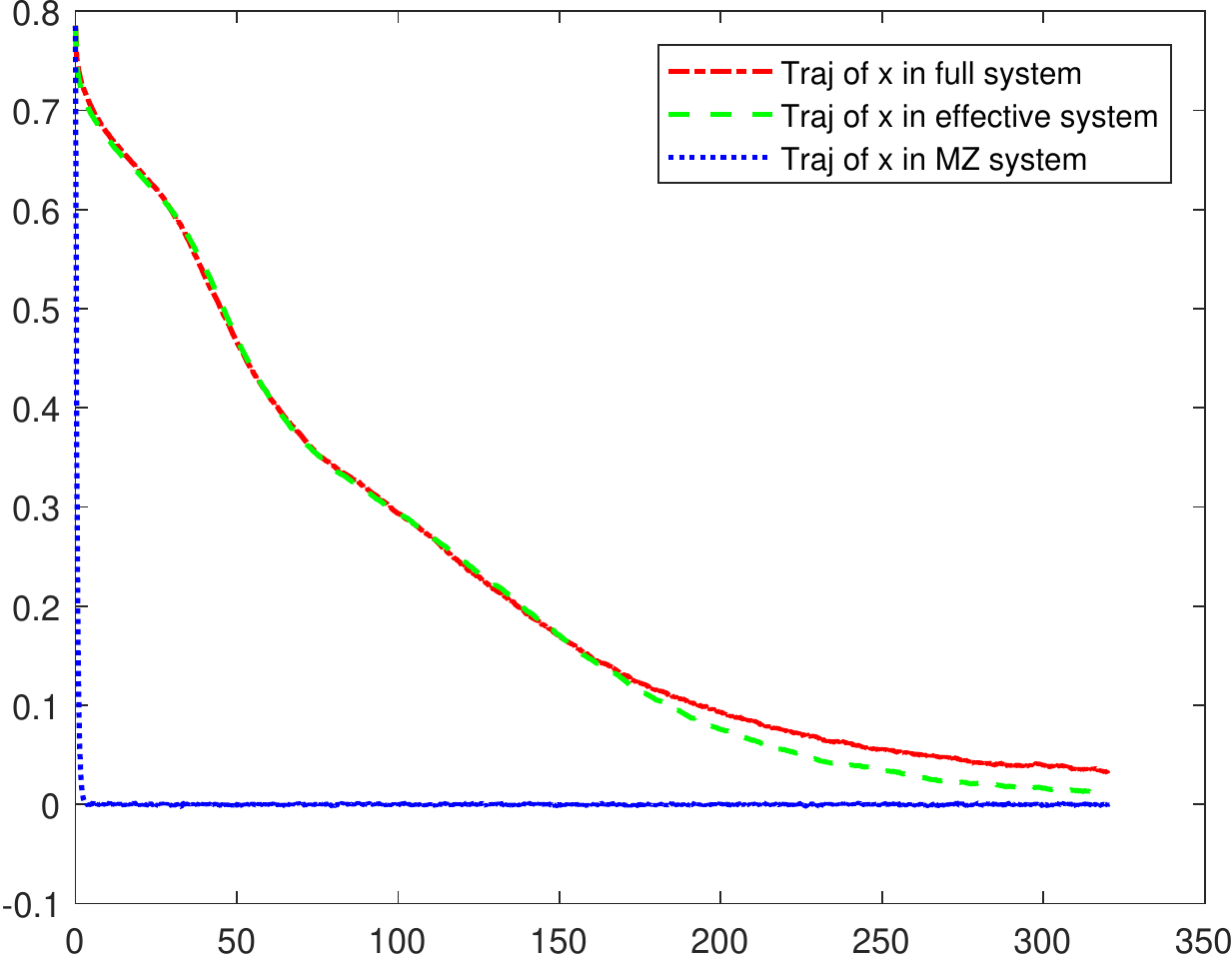}}
\vspace{- 0.2 pt}
\caption{Trajectories of $\mbf{h}_t$ averaged over 500 trajectories of
  \eqref{over_damped_Langevin}, \eqref{ApproxDynThermo} and \eqref{NoMemoryThermo}, with identical
  realisations of Brownian motion in each case. Parameters of $V$ are $\lambda=20$ and $\mu=2$, $\tau=2$, $\omega=10$, and the time step was $\Delta t=10^{-5}$, total time was $T=320$.\label{Fig:random_full_vs_revisedMZ}}
\end{figure}

\item\emph{Simulations with thermostat}.
  We next included the thermostat once more, considering sample averages of the effective dynamics
  \begin{align}
    d\mbf{h}_t&=-\frac{\mu\mbf{h}_t}{1+\tau^2\omega^2\cos^2(\omega \mbf{h}_t)}dt
                +\sqrt{\frac{2\beta^{-1}\Phi\Phi^T}{1+\tau^2\omega^2\cos^2(\omega \mbf{h}_t)}} d\mbf{B}_t\label{ApproxDynThermo}\\
    d\mbf{h}_t&=-\mu\mbf{h}_t+\sqrt{2\beta^{-1}\Phi\Phi^T} d\mbf{B}_t.\label{NoMemoryThermo}
  \end{align}
  For comparison, we again sampled the true dynamics of $\mbf{h}_t$ including the thermostat,
  i.e.
  \begin{equation}
    d\mbf{F}(\mbf{X}_t) = -\Phi \nabla V(\mbf{X}_t)  dt+\sqrt{2\beta^{-1}}\Phi d\mbf{B}_t.
    \tag{\ref{Liouville_eq2}}
  \end{equation}
  Simulations were performed at different choices of inverse temperature, and averages were taken
  over 500 identical realisations of the Brownian motion, aiming to minimise statistical error.
  In each case, initial conditions were chosen such that
  \[
    \cos(\omega X_0)=0,\qquad Y_0=\tau\sin(\omega Z_0).
  \]
  The results of this simulation are shown in Figure~\ref{Fig:random_full_vs_revisedMZ}.

  Consistent with our observation that the effective dynamics \eqref{ApproxDyn} is independent
  of $\beta$, we observe good agreements between \eqref{ApproxDynThermo} and the full dynamics
  for all choices of $\beta$, with increasing accuracy as $\beta\to\infty$. In comparison,
  \eqref{NoMemoryThermo} reflects the dynamical properties of the system poorly in all cases.
\end{enumerate}

\section{conclusion}
In this paper, we have employed the Mori--Zwanzig framework to rigorously derive an effective
equation for linear reaction coordinates, describing features of an underlying overdamped Langevin
dynamics.
Such models are appropriate for a variety of applications where we wish to capture
only limited aspects of a complex model, such as MD systems in the high friction limit.
The equation we derived enables us to understand the sources of error and thereby inform a choice
of effective dynamics which better captures dynamical features of the evolution which are not
well--represented by the dynamics of the effective potential alone.
We hope that this approach can serve to aid practitioners in understanding the sources of error in
a coarse--grained model, particularly in the presence of entropic barriers.

We validated our analytic results by considering a benchmark example of overdamped Langevin system
in a case where relaxation is impeded by a winding free energy barrier. This necessitated
the careful asymptotic treatment of interactions between resolved and unresolved variables in order to
correctly capture the dynamical behaviour of the system. In particular, although a time-scale
separation occurs within the system we considered, we nevertheless showed that
careful asymptotic analysis or dynamical sampling is required in general to ensure accuracy.
The approximate model we constructed based upon the equations we derived exhibited a drastic
improvement in predicting the dynamical behaviour of the reaction coordinates over the common
approach of using the effective potential alone to describe the dynamics.

Our work prompts several questions, which we hope to address in future:
\begin{enumerate}
\item \emph{Practical sampling algorithms and error analysis.} It would be of practical interest to devise an algorithm to generate
  effective dynamics based on the asymptotic approximation we considered here, and conduct a
  rigorous error analysis in this case.
\item \emph{Extensions to nonlinear coarse--grained variables and full Langevin dynamics.} In this
  work, we have considered the overdamped Langevin setting with linear reaction coordinates.
  It would be of significant interest to extend this analysis to nonlinear variables and
  full Langevin dynamics using our reliable asymptotic analysis approach in future.
\end{enumerate}

\appendix

\section{Proof of Theorem~\ref{thm:MZ_formulation}}
\label{sec:Proof}

In this section, we provide a proof of our main mathematical result, Theorem~\ref{thm:MZ_formulation}. The proof follows a similar strategy to other derivations using the Mori--Zwanzig formalism
given in the literature, notably \cite{Chorin2002a,GKH05,Hijon2010a}.

\subsection{Construction of `orthogonal' variables}
Our first step to construct variables which capture directions in phase space which are
`orthogonal' to those captured by the reaction coordinates $\mbf{F}$, i.e. a
foliation of the phase space.

Recall that $\Phi\in\mbb{R}^{m\times N}$ is a matrix of full rank by assumption, and therefore
it follows that the symmetric strictly positive
definite square root matrix $\Sigma:=\sqrt{\Phi\Phi^T}\in\mbb{R}^{m\times m}$ exists.
Recalling the construction of the singular value decomposition, we find that there exist
orthogonal matrices $U\in\mbb{R}^{m\times m}$ and $V\in\mbb{R}^{N\times N}$ such that
\begin{equation*}
  \Sigma^{-1}\Phi = UD V^T,\quad\text{where}\quad  D = \left(
    \begin{array}{cc}\mathrm{I}_m
      & 0
    \end{array}\right)\in\mbb{R}^{m\times N},
\end{equation*}
where $\mathrm{I}_m\in\mbb{R}^{m\times m}$ is the identity matrix, and $0$ denotes a submatrix of
zeros. Defining
\begin{equation*}
  E := \left(
    \begin{array}{cc} 0 & \mathrm{I}_{N-m}
    \end{array}\right)\in\mbb{R}^{(N-m)\times N},
\end{equation*}
where again $\mathrm{I}_{N-m}\in\mbb{R}^{(N-m)\times(N-m)}$ is an identity matrix and $0$ denotes a matrix of zeros, we set
\begin{equation*}
  \Psi := EV^T\in\mbb{R}^{(N-m)\times N},\qquad\text{and}\qquad\Phi^* :=\Phi^T\Sigma^{-2}\in\mbb{R}^{N\times m}
\end{equation*}
and it follows that
\begin{equation}\label{PartitionOfId}
  \Phi^*\Phi+\Psi^T\Psi = \mathrm{I}_{N\times N}.
\end{equation}

From the construction above, we see that $\Phi^*\Phi\in\mbb{R}^{N\times N}$
is an orthogonal projection acting on the phase space $\mbb{R}^N$, and the matrix $\Psi$ `selects' exactly the unresolved variables, so
that if $\Phi\mbf{x}=\mbf{h}$ and $\Psi\mbf{x} = \widetilde{\mbf{x}}$, we have
  \[
    \mbf{x}=\Phi^*\Phi \mbf{x}+\Psi^T\Psi\mbf{x}=\Phi^*\mbf{h}+\Psi^T\widetilde{\mbf{x}}.
  \]
  Given $\Phi$, we may use the construction of $\Psi$ in order to define the partition function
  $Z_\Phi:\mbb{R}^m\to\mbb{R}$,
  \begin{equation*}
    Z_\Phi(\mbf{h}):= \int  e^{-\beta V(\Phi^*\mbf{h}+\Psi^T\widetilde{\mbf{x}})} d\widetilde{\mbf{x}}.  \end{equation*}

  \subsection{Dyson--Duhamel principle}
  Next, we define
  \begin{equation*}
    \mbf{h}_t(\mbf{x}) = \mbb{E}\big[\mbf{F}(\mbf{x}_t)\,\big|\,\mbf{x}_0=\mbf{x}\big].
  \end{equation*}
  Applying the Feynman--Kac formula, we recall that $\mbf{h}_t$ solve the PDE
  \begin{equation*}
    \partial_t\mbf{h}_t = \calL \mbf{h}_t= -\nabla V\cdot\nabla\mbf{h}_t+\beta^{-1}\Delta\mbf{h}_t,
    \quad\text{with}\quad\mbf{h}_0(\mbf{x}) = \mbf{F}(\mbf{x}).
  \end{equation*}
  In semigroup notation, we will write $\mbf{h}_t = e^{t\calL}\mbf{F}$, and so the Feynman--Kac
  formula becomes
  \begin{equation}\label{FK1}
    \partial_t\mbf{h}_t = e^{t\calL}\calL\mbf{F}.
  \end{equation}
  Given mutually orthogonal projection operators $\calP$ and $\calQ$, applying the Dyson--Duhamel
  principle entails that we have the identity
  \begin{equation}\label{DysonDuhamel}
    e^{t\calL} = \int_0^t e^{(t-s)\calL}\calP\calL e^{s\calQ\calL} ds +e^{t\calQ\calL},
  \end{equation}
  which can be verified by differentiation with respect to $t$.
  Writing $\calL\mbf{F} = \calP\calL\mbf{F}+\calQ\calL\mbf{F}$ in \eqref{FK1}
  and applying \eqref{DysonDuhamel}, we find that
  \begin{equation}\label{FK2}
    \partial_t\mbf{h}_t = e^{t\calL}\calP\calL\mbf{F}
    +\int_0^t e^{(t-s)\calL}\calP\calL e^{s\calQ\calL}\calQ\calL\mbf{F}\,ds
    +e^{t\calQ\calL}\calQ\calL\mbf{F}.
  \end{equation}
  Our main focus will now be on the first two terms in \eqref{FK2}, since the latter term is
  $\calF_t$ as defined in \eqref{RandomForce}. To rewrite the former term,
  we apply the definition of projection operator $\calP$, the definitions of $\Phi$, $\Psi$ and
  $Z_\Phi$ and the chain rule, giving
    \begin{multline*}
\calP \calL \mbf{F}(\mbf{x})
=\Sigma^2\,\mbb{E}\big[ -\Sigma^{-2}\Phi\nabla V\,\big|\, \mbf{F}(\mbf{x})\big]\\
=\frac{\Sigma^2}{Z_\Phi\big(\mbf{F}(\mbf{x})\big)}\int -(\Phi^*)^T\nabla V(\Phi^*\,\mbf{F}(\mbf{x})+\Psi^T\widetilde{\mbf{x}}) e^{-\beta V(\Phi^*\,\mbf{F}(\mbf{x})+\Psi^T\widetilde{\mbf{x}}) } d\widetilde{\mbf{x}}
= -\Sigma^2\nabla\calS\big(\mbf{F}(\mbf{x})\big),
\end{multline*}
where we recall the definition of the effective potential $\calS$ given in \eqref{effective_potential}.
\subsection{Orthogonal dynamics}
Next, we consider the action of $e^{s\calQ\calL}$, which will subsequently allow us to treat the
integral term in \eqref{FK2}. We begin by noting that
\begin{equation*}
  \calQ \calL \mbf{F}(\mbf{x}) = \calL \mbf{F}(\mbf{x})-\calP \calL \mbf{F}(\mbf{x}) = -\Phi\nabla V(\mbf{x})+\Sigma^2\nabla\calS\big(\mbf{F}(\mbf{x})\big),
\end{equation*}
and define $\mbf{g}_t:\mbb{R}^N\to\mbb{R}^m$ to be the solution to
\begin{equation*}
  \partial_t \mbf{g}_t(\mbf{x}) = \calQ\calL\mbf{F}(\mbf{x})\cdot \nabla\mbf{g}_t(\mbf{x}),\qquad
  \mbf{g}_0(\mbf{x}) = \calQ\calL\mbf{F}(\mbf{x});
\end{equation*}
using semigroup notation, we write this as $\mbf{g}_s= e^{s\calQ\calL}\calQ\calL\mbf{F}$.

Under assumptions (1) and (2) made in Section~\ref{sec:prob_formulation}, it can be verified that
$\Sigma^2\nabla\calS\circ\mbf{F}$ is globally Lipschitz. Using this fact, it can therefore be shown
using the method of characteristics that $\mbf{g}_s$ exists and is a $C^1$ diffeomorphism on
$\mbb{R}^N$; for similar results in the Hamiltonian setting, see \cite{GKH05}.

Moreover, defining $\mbf{a}_s:=\mbf{g}_s(\mbf{x})\in\mbb{R}^m$ and $\mbf{b}_s:=\mbf{g}_s(\mbf{y})\in\mbb{R}^m$, then we use the fact that $\nabla\calS$ and $\nabla V$ are Lipschitz along with Young's
inequality to deduce that
\begin{align*}
  &\frac{d}{ds}\frac12\left|\mbf{a}_s-\mbf{b}_s\right|^2\\
  &\qquad=-\Big(\Phi\nabla V(\Phi^*\mbf{a}_s+\Psi^T\Psi\mbf{x})
    -\Phi\nabla V(\Phi^*\mbf{b}_s+\Psi^T\Psi\mbf{y})-\Sigma^2\nabla\calS(\mbf{a}_s)
    +\Sigma^2\nabla\calS(\mbf{b}_s)\Big)\cdot (\mbf{a}_s-\mbf{b}_s),\\
  &\qquad\lesssim \Big(\big|\Phi^*(\mbf{a}_s-\mbf{b}_s)+\Psi^T\Psi(\mbf{x}-\mbf{y})\big|+|\mbf{a}_s-\mbf{b}_s|\Big)|\mbf{a}_s-\mbf{b}_s|,\\
  &\qquad\lesssim |\mbf{a}_s-\mbf{b}_s|^2+|\mbf{x}-\mbf{y}|^2.
\end{align*}
Applying Gronwall's inequality in the usual way, it follows that there exists $\alpha>0$ such that
\begin{equation}\label{eq:NablaGBound}
  |\mbf{g}_s(\mbf{x})-\mbf{g}_s(\mbf{y})| \leq |\mbf{x}-\mbf{y}|\sqrt{s} e^{\alpha s},
  \qquad\text{and thus}\qquad |\nabla \mbf{g}_s(\mbf{x})|\leq \sqrt{s} e^{\alpha s}.
\end{equation}

\subsection{Memory integral}
Now that we have established properties of $\mbf{g}_s$, we return to the integral term in
\eqref{FK2}. For now, we fix $\mbf{x}\in\mbb{R}^N$, and set $\mbf{h}:=\mbf{F}(\mbf{x})$.

Consider $\calP\calL \mbf{g}_s$; using the partition of the identity constructed in
\eqref{PartitionOfId} we may write
\[
\begin{split}
  &\calP\calL \mbf{g}_s(\mbf{x})\\
  &\quad=
\frac{1}{Z_\Phi(\mbf{h})}
\int \Big(-\nabla \mbf{g}_{s}(\Phi^*\mbf{h}+\Psi^T\widetilde{\mbf{x}})
\cdot\nabla V(\Phi^*\mbf{h}+\Psi^T\widetilde{\mbf{x}})
+\beta^{-1}\Delta \mbf{g}_s(\Phi^*\mbf{h}+\Psi^T\widetilde{\mbf{x}})\Big) e^{-\beta V(\Phi^*\mbf{h}+\Psi^T\widetilde{\mbf{x}})} d\widetilde{\mbf{x}}\\
&\quad=
\frac{1}{Z_\Phi(\mbf{h})}
\int \Big(- \nabla \mbf{g}_{s}(\Phi^*\mbf{h}+\Psi^T\widetilde{\mbf{x}})\Phi^*\Phi\nabla V(\Phi^*\mbf{h}+\Psi^T\widetilde{\mbf{x}})+\beta^{-1}\Delta \mbf{g}_s(\Phi^*\mbf{h}+\Psi^T\widetilde{\mbf{x}})\Phi^*\Phi\\
&\hspace{40mm}- \nabla \mbf{g}_{s}(\Phi^*\mbf{h}+\Psi^T\widetilde{\mbf{x}})\Psi^T\Psi\nabla V(\Phi^*\mbf{h}+\Psi^T\widetilde{\mbf{x}})+\beta^{-1}\Delta \mbf{g}_s(\Phi^*\mbf{h}+\Psi^T\widetilde{\mbf{x}})\Psi^T\Psi\Big)\\
&\hspace{130mm}\times e^{-\beta V(\Phi^*\mbf{h}+\Psi^T\widetilde{\mbf{x}})} d\widetilde{\mbf{x}}.
\end{split}
\]
We collect terms involving matrix products with $\Phi$ and $\Psi$ separately, and using the
chain rule, we find that
\[
  \begin{split}
    \calP\calL \mbf{g}_s(\mbf{x})
&=\underbrace{\frac{1}{Z_\Phi(\mbf{h})}
\int \mathsf{div}_{\widetilde{\mbf{x}}}\left(\frac{1}{\beta}\nabla_{\widetilde{\mbf{x}}} \mbf{g}_s(\Phi^*\mbf{h}+\Psi^T\widetilde{\mbf{x}}) e^{-\beta V(\Phi^*\mbf{h}+\Psi^T\widetilde{\mbf{x}})}\right)d\widetilde{\mbf{x}}}_{=:T_1}\\
&\hspace{12mm}+\underbrace{\frac{1}{Z_\Phi(\mbf{h})}
\int \mathsf{div}_{\mbf{h}}\left(\frac{1}{\beta} \nabla_{\mbf{h}} \mbf{g}_s(\Phi^*\mbf{h}+\Psi^T\widetilde{\mbf{x}})\Sigma^2e^{-\beta V(\Phi^*\mbf{h}+\Psi^T\widetilde{\mbf{x}})}\right) d\widetilde{\mbf{x}}}_{=:T_2},
\end{split}
\]
where subscripts denote the variable with respect which derivatives are taken.
In particular, in the formula above,
$\mbf{g}_s(\Phi^*\mbf{h}+\Psi^T\mbf{x})$ is treated as a composition of functions. We now consider
each of the terms $T_1$ and $T_2$ separately.

To treat $T_1$, we apply the divergence theorem. Truncating the domain of integration to
$B_R(0)\subset\mbb{R}^{N-m}$, a ball of radius $R$ centred at $0$, and considering the limit as
$R\to\infty$, we have
\begin{equation}\label{T2_eq}
\begin{split}
T_1&=\frac{1}{Z_\Phi(\mbf{h})} \lim\limits_{R\rightarrow \infty}\int_{B_R(0) }
 \mathsf{div}_{\widetilde{\mbf{x}}} \left(\frac{1}{\beta} \nabla_{\widetilde{\mbf{x}}}\mbf{g}_s (\Phi^*\mbf{h}+\Psi^T\widetilde{\mbf{x}})e^{-\beta V(\Phi^*\mbf{h}+\Psi^T\widetilde{\mbf{x}})}\right)d\widetilde{\mbf{x}}\\
 &= \frac{1}{Z_\Phi(\mbf{h})}\lim\limits_{R\rightarrow \infty}
 \int_{\partial B_R(0) }
  \left(\frac{1}{\beta} \nabla_{\widetilde{\mbf{x}}}\mbf{g}_s (\Phi^*\mbf{h}+\Psi^T\widetilde{\mbf{x}}) e^{-\beta V(\Phi^*\mbf{h}+\Psi^T\widetilde{\mbf{x}})}\right)\cdot \vec{\nu}\,d\widetilde{S}.
\end{split}
\end{equation}
Applying \eqref{eq:NablaGBound} and the growth assumptions on $V$ to pass
to the limit, we see that $T_1=0$.

For $T_2$, we note that we may commute differentiation and integration, and so multiplying
and dividing by $Z_{\Phi}(\mbf{h})$, we obtain,
we obtain
\begin{equation}\label{T1_eq1}
  \begin{split}
T_2&=\frac{1}{Z_\Phi(\mbf{h})} \mathsf{div}_{\mbf{h}}
\left( \int \frac{1}{\beta} \nabla_{\mbf{h}} \mbf{g}_s(\Phi^*\mbf{h}+\Psi^T\widetilde{\mbf{x}})\Sigma^2e^{-\beta V(\Phi^*\mbf{h}+\Psi^T\widetilde{\mbf{x}})} d\widetilde{\mbf{x}}\right)\\
&= \frac{1}{Z_\Phi(\mbf{h})} \mathsf{div}_{\mbf{h}}
\left( \frac{Z_\Phi(\mbf{h})}{\beta}\frac{1}{ Z_\Phi(\mbf{h})} \int \nabla_{\mbf{h}} \mbf{g}_s(\Phi^*\mbf{h}+\Psi^T\widetilde{\mbf{x}})\Sigma^2e^{-\beta V(\Phi^*\mbf{h}+\Psi^T\widetilde{\mbf{x}})} d\widetilde{\mbf{x}}\right).
\end{split}
\end{equation}

\subsection{Memory kernel}
To complete our analysis, we must show the identity
\begin{equation}\label{KernelIdentity}
  \calM_s(\mbf{h}) = -\frac{1}{Z_\Phi(\mbf{h})} \int\nabla_{\mbf{h}} \mbf{g}_s(\Phi^*\mbf{h}+\Psi^T\widetilde{\mbf{x}})\Sigma^2e^{-\beta V(\Phi^*\mbf{h}+\Psi^T\widetilde{\mbf{x}})} d\widetilde{\mbf{x}},
\end{equation}
where we recall that $\calM_s$ was defined in \eqref{memory_eq2}.

Since we may again commute differentiation and integration, we use the product rule to write
\begin{equation*}
\begin{split}
&-\frac{1}{Z_\Phi(\mbf{h})}\int \nabla_{\mbf{h}} \mbf{g}_s(\Phi^*\mbf{h}+\Psi^T\widetilde{\mbf{x}})\Sigma^2e^{-\beta V(\Phi^*\mbf{h}+\Psi^T\widetilde{\mbf{x}})}
d\widetilde{\mbf{x}}\\
&\quad= -\frac{1}{Z_\Phi(\mbf{h})}\int \nabla_{\mbf{h}}\left( \mbf{g}_s(\Phi^*\mbf{h}+\Psi^T\widetilde{\mbf{x}})e^{-\beta V(\Phi^*\mbf{h}+\Psi^T\widetilde{\mbf{x}})}\right)\Sigma^2
d\widetilde{\mbf{x}}\\
 &\quad\qquad-\frac{1}{Z_\Phi(\mbf{h})}\int \mbf{g}_s(\Phi^*\mbf{h}+\Psi^T\widetilde{\mbf{x}})\otimes \beta\,(\Phi^*)^T\nabla V(\Phi^*\mbf{h}+\Psi^T\widetilde{\mbf{x}})e^{-\beta V(\Phi^*\mbf{h}+\Psi^T\widetilde{\mbf{x}})}\,\Sigma^2
 d\widetilde{\mbf{x}}
 \\
&\quad= -\frac{1}{Z_\Phi(\mbf{h})}
\nabla_{\mbf{h}}\left(\int \mbf{g}_s(\Phi^*\mbf{h}+\Psi^T\widetilde{\mbf{x}})e^{-\beta V(\Phi^*\mbf{h}+\Psi^T\widetilde{\mbf{x}})}d\widetilde{\mbf{x}}\right)\\
&\quad\qquad -\frac{1}{Z_\Phi(\mbf{h})}
\int\left( \mbf{g}_s(\Phi^*\mbf{h}+\Psi^T\widetilde{\mbf{x}}) \otimes \beta\,\Phi\nabla V(\Phi^*\mbf{h}+\Psi^T\widetilde{\mbf{x}})e^{-\beta V(\Phi^*\mbf{h}+\Psi^T\widetilde{\mbf{x}})}\right)d\widetilde{\mbf{x}},\\
&\quad=\underbrace{-\frac{1}{Z_\Phi(\mbf{h})}\nabla_{\mbf{h}}\Big(Z_\Phi(\mbf{h})\mbb{E}[\mbf{g}_s|\mbf{F}(\mbf{x})=\mbf{h}]\Big)}_{=:T_{11}}+\underbrace{\vphantom{\frac{1}{Z_\Phi(\mbf{h})}}\beta\,\mbb{E}[\mbf{g}_s\otimes\calL\mbf{F}|\mbf{F}(\mbf{x})=\mbf{h}]}_{=:T_{12}}.
\end{split}
\end{equation*}
Next, we recall that $\mbf{g}_s = e^{s\calQ\calL}\mbf{F}$, and $\calP\calQ=0$, so
\begin{equation}\label{eq:PesQL=0}
  \mbb{E}\big[\mbf{g}_s\big|\mbf{F}(\mbf{x})\big]=\mbb{E}\big[e^{s\calQ\calL}\calQ\calL\mbf{F}\big|\mbf{F}(\mbf{x})\big]=\left(\calP\calQ\calL e^{s\calQ\calL}\mbf{F}\right)(\mbf{x})=0,
\end{equation}
and therefore $T_{11}=0$. To treat $T_{12}$, we note that since $\calP+\calQ=\calI$, we may
split $T_{12}$ into
\begin{equation*}
\begin{split}
  T_{12} &= \beta\,\mbb{E}[e^{s\calQ\calL}\calQ\calL\mbf{F}\otimes\calQ\calL\mbf{F}|\mbf{F}(\mbf{x})=\mbf{h}]
  +\beta\,\mbb{E}[e^{s\calQ\calL}\calQ\calL\mbf{F}\otimes\calP\calL\mbf{F}|\mbf{F}(\mbf{x})=\mbf{h}]\\
  &= \calM_s(\mbf{h})
  +\beta\,\calP e^{s\calQ\calL}\calQ\calL\mbf{F}\otimes\calP\calL\mbf{F}.
\end{split}
\end{equation*}
Once again, the latter term vanishes thanks to \eqref{eq:PesQL=0}, and so we have proved identity
\eqref{KernelIdentity}.

\subsection{Conclusion of the proof}
Applying identity \eqref{KernelIdentity} to \eqref{T1_eq1} and using the product rule and the
definition of the effective potential given in \eqref{effective_potential}, we find that
\begin{equation*}
  T_2 = \frac{1}{Z_\Phi(\mbf{h})} \mathsf{div}
  \left( -\frac{Z_\Phi(\mbf{h})}{\beta}\calM_s(\mbf{h})\right) = \calM_s(\mbf{h})\nabla\calS(\mbf{h})
  -\frac{1}{\beta}\mathsf{div}\calM_s(\mbf{h}).
\end{equation*}
Combining our analysis of each of the terms, we have therefore shown that
\begin{equation*}
\partial_t\mbf{h}_t = -\Sigma^2\nabla \calS(\mbf{h}_t)+\int_0^t\calM_s(\mbf{h}_{t-s})\nabla\calS(\mbf{h}_{t-s})-\frac{1}{\beta}\mathsf{div}\calM_s(\mbf{h}_{t-s})\,ds +\calF_t
\end{equation*}
Hence, we prove the theorem.

\section*{Acknowledgment}
We would like to thank Dr Xiantao Li for helpful suggestions and encouragement during this project.
\bibliographystyle{abbrv}
\bibliography{MZbenchmark}
\end{document}